\documentclass[lettersize,journal]{IEEEtran}
\usepackage{amsmath,amssymb}
\usepackage{graphicx} 
\usepackage{amsthm}
\usepackage{mathtools}
\usepackage{bm}
\usepackage{bbm} 
\usepackage{cite}
\usepackage{graphicx}
\usepackage{booktabs}
\usepackage{multirow}
\usepackage{float}
\usepackage[table]{xcolor} 
\usepackage[labelfont=small]{caption} 
\usepackage{subcaption}

\usepackage{enumerate}

\usepackage{tcolorbox}
\tcbset{
    mygraybox/.style={
        colback=gray!10, 
        colframe=gray!80, 
        sharp corners, 
        boxrule=1pt, 
        width=\linewidth, 
        enlarge left by=-1mm, 
        enlarge right by=-1mm, 
    }
}

\usepackage{color,soul}
\usepackage{svg}

\newcommand{\set}[1]{\mathcal{#1}}
\newcommand{\norm}[1]{\left\lVert#1\right\rVert}

\usepackage[noabbrev,capitalize]{cleveref}

\crefname{equation}{}{}
\Crefname{equation}{}{}
\crefname{section}{}{}
\Crefname{section}{}{}

\DeclareMathOperator{\diag}{diag}

\theoremstyle{definition} 

\newtheorem{assumption}{Assumption}

\theoremstyle{plain} 

\newtheorem{definition}{Definition}

\theoremstyle{remark} 
\newtheorem{remark}{Remark}

\title{
Analysis of Data Value in Stochastic Optimal Power Flow for Distribution Systems
}

\author{Mehrnoush Ghazanfariharandi and Robert Mieth}

\usepackage{amsmath,amsfonts}
\usepackage{algorithmic}
\usepackage{array}
\usepackage[caption=false,font=normalsize,labelfont=sf,textfont=sf]{subfig}
\usepackage{textcomp}
\usepackage{stfloats}
\usepackage{url}
\usepackage{verbatim}
\usepackage{graphicx}
\usepackage[font=small]{caption}
\def\BibTeX{{\rm B\kern-.05em{\sc i\kern-.025em b}\kern-.08em
    T\kern-.1667em\lower.7ex\hbox{E}\kern-.125emX}}

\usepackage{pict2e}
\usepackage{keyval}
\usepackage{calc}
\usepackage{fp}
\usepackage{diagbox}
\usepackage{booktabs}
\usepackage{tabularx,colortbl, makecell, caption, multirow}
\usepackage{boldline} 
\usepackage{tabu,multirow}

\usepackage{placeins}
\newcommand{\subparagraph}{}
\usepackage{titlesec}

\begin{document}

\bstctlcite{IEEE:BSTcontrol}

\renewcommand{\thesection}{\Roman{section}} 
\renewcommand{\thesubsection}{\Alph{subsection}}

\maketitle
\thispagestyle{empty}
\pagestyle{empty}


\begin{abstract}
\noindent 
The rise of advanced data technologies in electric power distribution systems enables operators to optimize operations but raises concerns about data security and consumer privacy. Resulting data protection mechanisms that alter or obfuscate datasets may invalidate the efficacy of data-driven decision-support tools and impact the value of these datasets to the decision-maker.  
This paper derives tools for distribution system operators to enrich data-driven operative decisions with information on data quality and, simultaneously, assess  data usefulness in the context of this decision. 
To this end, we derive an AC optimal power flow model for radial distribution systems with data-informed stochastic parameters that internalize a data quality metric. 
We derive a tractable reformulation and discuss the marginal sensitivity of the optimal solution as a proxy for data value. 
Our model can capture clustered data provision, e.g., from resource aggregators, and internalize individual data quality information from each data provider.
We use the IEEE 33-bus test system, examining scenarios with varying photovoltaic penetration and load scenarios, to demonstrate the application of our approach and discuss the relationship between data quality and its value. 
\end{abstract}

\section{Introduction}
\IEEEPARstart{E}{nhanced} visibility of electrical distribution grids enabled by modern data collection and communication technology facilitates sophisticated operational strategies for distribution system operators (DSOs). 
Specifically, the availability of high-resolution voltage and current data in conjunction with an increased deployment of controllable distributed resources enable DSOs to guide system operations towards achieving safe and efficient optimal power flow (OPF) solutions \cite{dall2017chance}. 
However, the collection of this data can compromise the privacy of electricity consumers \cite{giaconi2021smart}.
Resulting privacy concerns towards digital data collection and smart meters can pose significant adoption barriers for data sharing technology in distribution systems and, as a result, delay the integration of renewable and distributed energy resources \cite{erkin2013privacy,le2020ethical}.
The ongoing push to utilize data for system operations while protecting individual data privacy and data integrity has motivated the development of data obfuscation methods that balance data usefulness and privacy protection \cite{giaconi2021smart,currie2023data}.
However, dealing with obfuscated data and assessing its usefulness for DSOs is complicated by the complexity of decision-making processes influenced by power flow physics.
As a result, tools for quantifying data value and informing investments in data acquisition are lacking. 
Motivated by this gap, this paper presents a direct data-valuation approach for data-driven active distribution system management. 
Building on a data quality metric introduced in \cite{mieth2023data}, 
we derive a chance-constrained AC OPF formulation that internalizes information on data quality from various sources (e.g., aggregators) to optimally manage distributed resources and to obtain situation- and location-aware data value proxies. 

Prompted by the need to deal with increasing
levels of uncertainty, mainly from uncontrollable high-wattage grid-edge resources such as rooftop solar, residential battery systems, and EV charging, various stochastic and data-driven decision-support tools have been proposed for distribution systems  \cite{dall2017chance,guo2018data,hassan2020stochastic,mieth2018data} and power systems in general \cite{roald2023power}.
Most of these approaches require data; Historical or real-time observations and measurements that can be used to estimate relevant parameters, train models, or quantify the statistics of random processes, e.g., related to load and renewable energy forecasts. 
Data availability and \textit{data quality} directly impact the efficacy of these methods. This fact is well-known in the context of state-estimation \cite{geth2023data} and has been studied to improve noise representation \cite{vanin2023exact} and the effect of uncertain renewable power injection \cite{zhang2019interval}.
Data-driven stochastic optimal power flow models for the optimal dispatch of distributed resources, on the other hand, which, e.g.,
rely on data to estimate the distribution of the underlying uncertainty \cite{mieth2018data,hassan2020stochastic} or directly pursue a data-driven approach that uses observed or sampled data points \cite{dall2017chance,guo2018data}, 
are typically ignorant to the quality of the available data. 
As a result, the role of the data as a central input to the problem remains overlooked. 

Most proposals for data-driven distribution system operations, e.g., as in  \cite{dall2017chance,guo2018data,hassan2020stochastic,mieth2018data}, assume that the necessary data streams are readily available at no cost to the decision-maker, e.g., the DSO.
Yet, the collection and preparation, or acquisition, of data at sufficient quality incurs expenses, either directly from data storage, processing, and transmission \cite{ren2018datum}, or indirectly from \textit{privacy loss} \cite{bessa2018data}.
As a result, the DSO must consider compensation or incentive payments to data owners to gain access to relevant data streams \cite{han2021monetizing}. If data owners can quantify their cost of sharing data, auction mechanisms and optimal data allocation algorithms, e.g., based on Shapely Value computations \cite{koutsopoulos2015auctioning,agarwal2019marketplace}, have been shown to deliver interpretable economic insights  on how data should be valued from the perspective of the data user. 
Such market-based approaches to data collection and sharing have also been proposed and discussed for power system applications \cite{han2021monetizing,goncalves2020towards,pinson2022regression,yassine2015smart,acharya2022false}. 
However, they require a bidding process managed by a third-party entity and an assessment of data cost by the data provider that informs their bidding participation.

In the context of distribution system operations, typical data owners are electricity consumers that are often required to share data through contractual agreements, e.g., with their metering service provider or an aggregator of distributed resources.
However, because this data encodes sensitive private information \cite{molina2010private}, many local regulators enforce prudent data protection laws that restrict or completely prohibit the general use or sharing of this data \cite{EU-GDPR,DOE2015voluntary}.
To unlock the benefits of data-driven operational processes, data protection mechanisms based on algorithmic obfuscation \cite{giaconi2021smart} or differential privacy \cite{toubeau2022privacy,dvorkin2020differentially} have been studied and can be expected to become a central enabler to emerging digitalized data-driven DSOs \cite{currie2023data}.

Naturally, data obfuscation such as differential privacy alters the data and impacts any downstream decision processes that use this data. 
In addition, data alteration can be intentional to ensure data privacy, as discussed above, or unintentional as a result of insufficient or noisy measurement infrastructure \cite{vanin2023exact,geth2023data}. 
In this case, information loss occurs before data is shared and is independent of how data is transmitted between the data provider and the data user. 

Motivated by the need for a better understanding of data quality and data value in distribution systems, we derive a chance-constrained AC OPF formulation for distribution system operations with uncertain load and renewable availability that uses data to estimate the distribution of the relevant uncertain parameters. 
Our proposed model has the following properties and makes the following contributions:
\begin{enumerate}[(i)]
    \item  It allows the decision-maker to internalize information on potential data obfuscation or alteration relative to the original data source, i.e., our model internalizes information on \textit{data quality}.
    \item It imputes a proxy for data value depending on data quality. The resulting data value provides interpretable insights for the DSO on the usefulness or importance of data \textit{in the context} of the decision-making problem.
    \item It generalizes the model in \cite{mieth2023data} to allow a single data provider to submit data on various features throughout the system. This enables the model to better reflect real-world grid operations that typically involve aggregators or similar third-party service providers.
\end{enumerate}
Our formulation utilizes results from \cite{mieth2023data}, which first presented the idea of internalizing heterogeneous data quality information into a decision-making problem via distributionally robust optimization (DRO). 
Our paper presents the required modifications for an AC OPF model (the work in \cite{mieth2023data} only discussed a DC OPF model with a linear recourse policy) and generalizes the results from \cite{mieth2023data} to accommodate datasets that provide information on a collection of uncertain parameters, e.g., collected by an aggregator.
As a result, this paper enables an additional layer of flexibility in data sourcing and utilization.
We demonstrate our approach with numerical case studies for various PV deployment and load scenarios.

A central motivation of this paper, is to analyze the usefulness or value of a dataset (i) as a function of its quality (i.e., how well the data reflects the true parameter) and (ii) in the context of the decision-making problem.
Unlike existing methods that require ex-post data valuation strategies, such as the Shapley Value \cite{han2021monetizing, goncalves2020towards,pinson2022regression}, our approach directly computes the marginal value of data quality, eliminating the need for additional ex-post calculations.
Moreover, data valuation based on combinatorial approaches like Shapely Value is computationally expensive to compute or estimate, or is sensitive to how the data is divided into training and testing. 
The method discussed in this paper offers a direct approach that emphasizes data quality, especially in contexts involving privacy protection and noisy data. 

The remainder of this paper is structured as follows: Section~\ref{sec:Data-driven system operations}  describes how we set up and solve the data-driven system operations. Section~\ref{sec:Internalizing data quality} explores the internalization of data quality, presenting the methodologies used to incorporate quality metrics into operational strategies. In Section~\ref{sec:casestudy}, we demonstrate our model on the IEEE 33-bus test system.  Section~\ref{sec:Conclusion} concludes the paper.

\section{Data-driven system operations}
We first formulate a data-driven AC OPF model for distribution system operations under uncertainty along the lines of \cite{dall2017chance}. We refer to Table~\ref{tab:symboldescription} for an overview of the main symbols and notations. 

\label{sec:Data-driven system operations}
\subsection{Problem formulation}
\label{ssec:problem_formulation}
We take the perspective of a DSO who solves an OPF problem for a distribution network with $N$ nodes indexed by $n = 1,..., N$. 
For given active and reactive load $\bm{p}_l=[p_{l,1},...,p_{l,N}]$, $\bm{q}_l=[q_{l,1},...,q_{l,N}]$ and maximum available RES injections $\bm{p}_{av}=[p_{av,1},...,p_{av,N}]$, the OPF computes optimal setpoints $\bm{\alpha}=[\alpha_1,...,\alpha_n]$, and $\bm{q}_c=[q_{c,1},...,q_{c,N}]$ of RES inverters and schedules controllable DERs $\mathbf{p}_B=[p_{B,1},...,p_{B,N}]$ and $\mathbf{q}_B=[q_{B,1},...,q_{B,N}]$.
For each inverter at bus $n$, $\alpha_n\in[0,1]$ controls the active power output and denotes the fraction of curtailed available power $p_{av,n}$. The same inverter can supply reactive power $q_{c,n}$ limited by its apparent power rating $S_n$ such that:
\begin{equation}
    ((1-\alpha_n)p_{av,n})^2 + q_{c,n}^2 \le S_n^2.
\label{eq:inverter_constraint}
\end{equation}

The objective of the DSO is to minimize a cost function $C$ (e.g., minimizing energy payments or minimizing power losses), while ensuring compliance with system voltage limits $V_{max},V_{min}$.
We collect the nodal voltage magnitudes in vector $\bm{\rho}:= [v_1,...,v_N] \in \mathbb{R}^N$ and leverage an approximate linear relationship between voltage magnitudes and active and reactive nodal power injections as:
\begin{equation}
\bm{\rho} \approx \mathbf{R}\mathbf{p} + \mathbf{B}\mathbf{q} + \mathbf{a}.
\label{eq:linear_voltage_magnitudes}
\end{equation}
There exist various effective methods in the literature to compute the parameters of $\mathbf{R}$, $\mathbf{B}$, and $\mathbf{a}$ \cite{dall2017chance}.
For this paper,
we rely on the established \textit{LinDistFlow} model \cite{baran1989optimal} for radial networks.
In this power flow model $\bm{\rho}$, in fact, represents the \textit{square} of the voltage magnitudes. However, we note that the method we derive below is independent of how $\mathbf{R}$, $\mathbf{B}$, and $\mathbf{a}$ are defined. Accommodating $\bm{\rho}$ as the squared voltage only changes the definition of the voltage limits.
We use \cref{eq:linear_voltage_magnitudes} to define the following vector-valued function:
\begin{align}
g_\rho(\bm{\alpha},\bm{q}_c,\bm{p}_B,\bm{q}_B,\bm{\delta}) &= \mathbf{R}\big((\bm{I}-{\rm diag}(\bm{\alpha}))\mathbf{p}_{av} - \mathbf{p}_l + \mathbf{p}_B \big)\nonumber \\& +\mathbf{B}\big(\mathbf{q}_c - \mathbf{q}_l +\mathbf{q}_B)\big) +\mathbf{a},
\end{align}
which computes the voltage magnitudes from variables $\bm{\alpha},\bm{q}_c,\bm{p}_B, \bm{q}_B $ and parameters $\mathbf{p}_{av}, \mathbf{p}_l , \mathbf{q}_l$, which we collect in vector $\bm{\delta}\coloneqq(\mathbf{p}_{av}, \mathbf{p}_l , \mathbf{q}_l)$.

At the time the DSO solves the OPF problem, the exact values of $\bm{\delta}$ are unknown,
e.g., due to forecast uncertainty \cite{dall2017chance}, imperfect measurements \cite{geth2023data}, or state-estimation with imperfect data \cite{vanin2023exact}.  
We highlight that in this paper we focus on uncertainty from forecasted load and renewable availability. 
Also, other parameters, such as line parameters or topology information, may be subject to uncertainty in practice but are considered known here. 

To capture this uncertainty, we model
$\mathbf{p}_{av}$ , $\mathbf{p}_l$ and $\mathbf{q}_l$ as random variables and the DSO must solve the OPF problem as a stochastic program with probabilistic performance guarantees.
Several such approaches have been proposed in the recent literature, e.g., \cite{dall2017chance,mieth2018data,li2018distribution}. In this paper we build our discussion on results from \cite{dall2017chance}, modified to enforce a joint chance constraint on the voltage limits:
\allowdisplaybreaks
\begin{subequations}
\label{eq:main}
\begin{flalign}
&\min_{{\bm{\alpha}, \bm{q}_c,\bm {p}_B,\bm{q}_B,\bm{\delta}}}\quad\mathbb{E}_\mathbbmss{Q}\left[C\left(\bm{\alpha},\bm{q}_c,\bm{p}_B,\bm{q}_B,\bm{\delta}\right)\right] \label{eq:DRO_OBJ1} \\
&\text{s.t.} \nonumber \\
&\hspace{0em}
 \mathrm{Pr}_{\mathbb{Q}}\! \left\{\begin{aligned}
 \!\!\!g_{\rho,n} (\!\bm \alpha,\bm q_c,\bm p_B,\bm{q}_B,\bm \delta\!) \!\le \!V_{max}\!\! \quad\!  n \! = \!1,..., N\\
 \!\!-g_{\rho,n} (\!\bm \alpha,\bm q_c,\bm p_B,\bm{q}_B,\bm \delta\!) \!\!\le \!\! -V_{min}\!\! \quad \!  n \! =\! 1,..., N\label{eq:jointccfirst}\!
\end{aligned}\right\} \!\!\ge \!\! 1\!-\eta^{vol}\\
&\hspace{0em}   \mathrm{Pr}_{\mathbb{Q}}[ ((1\!\!-\!\alpha_n\!) p_{av,n})^2\!\!+ q_{c,n}^2\!-\!S_n^{2}\! \le \!0 ] \!\ge 1-\eta^{inv} \label{eq:base_invertor}\quad n \! = \!1,...,N \\
&\hspace{0em} 0\le \alpha_n \le 1 \qquad n = 1,...,N  \label{eq:deterministiccon3} \\
&\hspace{0em} p_{B,n}^{min}\le p_{B,n} \le p_{B,n}^{max} \qquad n = 1,...,N \label{eq:deterministiccon1} \\
&\hspace{0em} q_{B,n}^{min}\le q_{B,n} \le q_{B,n}^{max} \qquad n = 1,...,N \label{eq:deterministiccon5} 
\end{flalign}%
\label{prob:main_op}%
\end{subequations}%
\allowdisplaybreaks[0]%
Objective \cref{eq:DRO_OBJ1} minimizes the expectation of the cost function for a given distribution $\mathbb{Q}$ of $\bm{\delta}$. 
Constraints \eqref{eq:jointccfirst} and \eqref{eq:base_invertor} ensure that for a given decision of inverter and DER setpoints $\bm \alpha,\bm q_c,\bm p_B, \bm q_B $ voltage constraints, and inverter capacity limits are met with a probability of $1-\eta^{vol}$ and $1-\eta^{inv}$, respectively.
Here, $g_{\rho,n}(\cdot)$ denotes the n-th element of $g_{\rho}(\cdot)$. 
Constraint \cref{eq:deterministiccon3} limits the share of curtailed renewable injection to $[0,1]$ and constraint \cref{eq:deterministiccon1} and \cref{eq:deterministiccon5} enforce technical DERs constraints. 

In this paper we define the cost function as the DSO's cost:
\begin{equation}
\begin{aligned}
C(\bm{\alpha},\bm{q}_c,\bm{p}_B,\bm{q}_B,\bm{\delta})&\! =\!\! 
\sum_{n=1}^{N}\!\big( c_n[p_{l,n} - \!p_{B,n} \!- (\!1-\alpha_n) p_{av,n}]^+\\& +  d_n [(1-\alpha_n) p_{av,n} - p_{l,n} + p_{B,n}]^+ \\&+ e_n ||q_{c,n}|+|q_{B,n}|| +  h_n \alpha_n p_{av,n}\big).
\end{aligned}
\end{equation}
Parameters $\bm c, \bm d, \bm e$, and $\bm h$ capture the cost of buying power from the grid, DSO feed-in tariff payments, payments for reactive power supply  from
inverters and DERs, and reimbursements for active power curtailment, respectively.

\begin{table}
\centering
\caption{Nomenclature}
\label{tab:symboldescription}
\begin{tabular}{c p{5.5cm}}
\hline
\textbf{Symbol} & \textbf{Description} \\
\hline
$c_n$ & Cost of buying power from the grid at node $n$\\
$d_n$ & DSO feed-in tariff payments at node $n$ \\
$e_n$ &   Payments for reactive
power supply at node $n$\\
$F$ & Number of data providers (data clusters) \\
$h_n$ &  Reimbursements for
active power curtailment at node $n$\\
$I$ & Number of samples \\
$K$& Number of voltage constraints $K=2N$ \\
$N$ & Number of nodes \\
$N_f$ & Number of nodes in cluster $f$ \\
$p_{av,n}$ & Maximum available RES active power at node $n$ \\
$p_{B,n}$ & Active power (dis)charging of DER at node $n$ \\
$p_{l,n}$ & Active power load at node $n$\\
$q_{B,n}$ & Reactive power (dis)charging of DER at node $n$ \\
$q_{c,n}$ & Reactive power provided by RES inverter at node $n$ \\
$q_{l,n}$ & Reactive power load at node $n$\\
$S_{n}$ & Rated apparent power of RES inverter at node $n$ \\
$v_n$ & Voltage magnitude at node $n$  \\
$V_{min}, V_{max}$ & Lower and upper limits for voltage magnitudes \\
$\alpha_{n}$ & Fraction of RES active power curtailed at node $n$ \\
$\bm{p}_{av}$ & Vector of maximum available RES injections\\
$\bm{p}_l$ & Vector of active power load\\
$\bm{q}_l$ & Vector of reactive power load\\
$\check{\bm{p}}_{av}$ & Vector of maximum available RES injections forecast \\
$\check{\bm{p}}_{l}$ & Vector of active power load forecast \\
$\bm{q_{c}}$ &  Vector of reactive power provided by RES  \\
$\check{\bm{q}}_{l}$ & Vector of reactive power load forecast \\
$\bm{x}$ & Vector collecting model decision variables\\
$\bm{\alpha}$ & Vector of fraction of active power curtailed by RES \\
$\bm{\delta}$ & Vector of uncertain parameters\\& $\bm{\delta}\coloneqq(\mathbf{p}_{av}, \mathbf{p}_l , \mathbf{q}_l)$ \\
$\bm{\rho}$ & Vector of nodal voltage magnitudes  \\
$\epsilon_f$ & Data quality of cluster $f$\\
$\eta^{inv}$ & Risk-level for inverter constraint violation\\
$\eta^{vol}$ & Risk-level for voltage constraint violation\\
$\varphi^{inv}_n$ & Auxiliary variable for CVaR related to RES capacity at node $n$\\
$\varphi^{vol}$ &Auxiliary variable for CVaR related to voltages\\
$\varpi_n^{inv}$ & Auxiliary variable related to RES capacity at node $n$\\
$\varpi^{vol}$ & Auxiliary variable related to voltages\\
$[\cdot]^+$ & $\max\{\cdot,0\}$\\
\hline
\end{tabular}
\end{table}

\subsection{Data-driven solution}
\label{sec:data-driven_solution}
The problem in \cref{prob:main_op} can be solved in a data-driven manner using historically recorded or sampled realizations $\{\widehat{\bm{\delta}}_i\}_{i=1}^{I}$ of $\bm{\delta}$.
With access to such samples, the DSO can replace objective \cref{eq:DRO_OBJ1} with its sample average approximation 
\begin{equation}
\!\mathbb{E}_{\widehat{\mathbbmss{P}}_I}\!\big[C\left(\bm{\alpha},\bm{q}_c,\bm{p}_B,\bm{q}_B,\bm{\delta}\right)\big] \!\!= \frac{1}{I}\sum_{i=1}^I C(\bm{\alpha},\bm{q}_c,\bm{p}_B,\bm{q}_B,\widehat{\bm{\delta}}_i),
\label{eq:objective_sample}
\end{equation}
where $\widehat{\mathbbmss{P}}_I$ denotes the empirical distribution of the $I$ available samples.
Similarly, chance constraints \cref{eq:jointccfirst,eq:base_invertor} allow a tractable data-driven reformulation using conditional value-at-risk (CVaR).
For \cref{eq:jointccfirst} this reformulation takes the form 
\begin{equation}
    \frac{1}{I}\sum_{i=1}^I\big[\max _{k = 1,...,K}[\langle \bm{a}_k, \bm{\delta}\rangle + c_k ] + \varphi^{vol}\big]^+\! \!\le \!\varphi^{vol} \eta^{vol} \quad  i=1,...,I, 
\label{eq:joint_cc_sample}
\end{equation}
where $\bm{a}_k$ and $c_k$ collect decision variables related to the voltage constraints and $K=2N$, and $\varphi^{vol}$ is an auxiliary decision variable.
We refer to Appendix~\ref{appendix:a} for more details on the derivation in the context of this paper and to \cite{dall2017chance,mieth2023data} for further details on the CVaR reformulation. 
Similarly, \cref{eq:base_invertor} can be written as
\begin{flalign}
   &  \frac{1}{I}\sum_{i=1}^I\big[
    ((1-\alpha_n) \widehat{p}_{av,n,i})^2 \!+ \! (q_{c,n})^2 \!- \!S_{n}^{2}  + \varphi^{inv}_n
    \big]^+ \!\!\!\le \! \varphi^{inv}_n \eta^{inv}\nonumber \\ &
        \! i=1,...,I,\quad n=1,...,N,
\label{eq:inverter_sample}
\end{flalign}
 where $\varphi^{inv}_n\in\mathbb{R}\ n=1,...,N$ are auxiliary decision variables.

The quality of the data-driven decision using reformulations \cref{eq:objective_sample,eq:joint_cc_sample,eq:inverter_sample} of \cref{eq:DRO_OBJ1,},\cref{eq:jointccfirst}, and \cref{eq:base_invertor}, respectively, relies on the ability of dataset $\{\widehat{\bm{\delta}}_i\}_{i=1}^I$ to accurately capture the true distribution of the uncertain $\bm{\delta}$. 
Therefore, $\{\widehat{\bm{\delta}}_i\}_{i=1}^I$ in itself is a crucial input to the decision-making process. 
In the following section, we present a modification of the data-driven solution of \cref{prob:main_op} that allows the system operator to evaluate the sensitivity of the solution to the input data $\{\widehat{\bm{\delta}}_i\}_{i=1}^I$ and model data accuracy obtained from various data sources.

\section{Internalizing data quality}
\label{sec:Internalizing data quality}
Assume that the DSO obtains the required input data $\{\widehat{\bm{\delta}}_i\}_{i=1}^I$ from multiple data providers, for example, metering service providers, aggregators, forecasting services, or smart home services.
It can be expected that each data provider follows privacy protection rules and subjects their data to some form of alteration or obfuscation, e.g., using differential privacy \cite{mieth2023data,giaconi2021smart,currie2023data}. 
As a result, the data available to the DSO might be noisy or biased and, as a result, may not accurately represent the true underlying distribution of the uncertain parameter. 
The same effect can occur if data measurements are biased, scarce, or noisy, regardless of any privacy protection mechanism. 
We highlight, that this process is independent of the actual data \textit{transmission} between data provider and data user, which must also be private and secure but can be assumed to be lossless for practical purposes.

The data provider can quantify the resulting discrepancy between the distribution supported by the submitted data and the true underlying distribution using a distributional distance metric.
The Wasserstein metric is a pertinent choice because of its properties for data-driven problems and in optimization \cite{mieth2023data}.
We therefore assume that each data provider computes or estimates the Wasserstein distance between the empirical distribution supported by their data and the true distribution of the parameter of interest.
By submitting this information alongside the data, each data provider can signal the quality of their data. 
A large Wasserstein distance indicates a higher level of data alteration, i.e., a lower data quality, and vice versa. 
If the data provider can compare the unaltered and altered datasets, then the Wasserstein distance between the empirical distributions supported by these datasets can be computed exactly. 
Moreover, data protection mechanisms that add noise to the data, like differential privacy, allow a direct computation of the upper bound of the distributional shift resulting from the added noise \cite{mieth2023data}.
If the true distribution cannot be observed, e.g., in the case of noisy or inadequate measurements, data quality can be estimated as a confidence upper bound of the Wasserstein distance using statistical methods as discussed in \cite{mohajerin2018data,le2024universal}.
We also refer to \cite{mieth2023data} for additional discussion.
The following section shows how this data quality information can be internalized in the data-driven OPF. 

\subsection{Multi-source distributional ambiguity}
\begin{figure}
    \centering
    \includegraphics[width=1\linewidth]{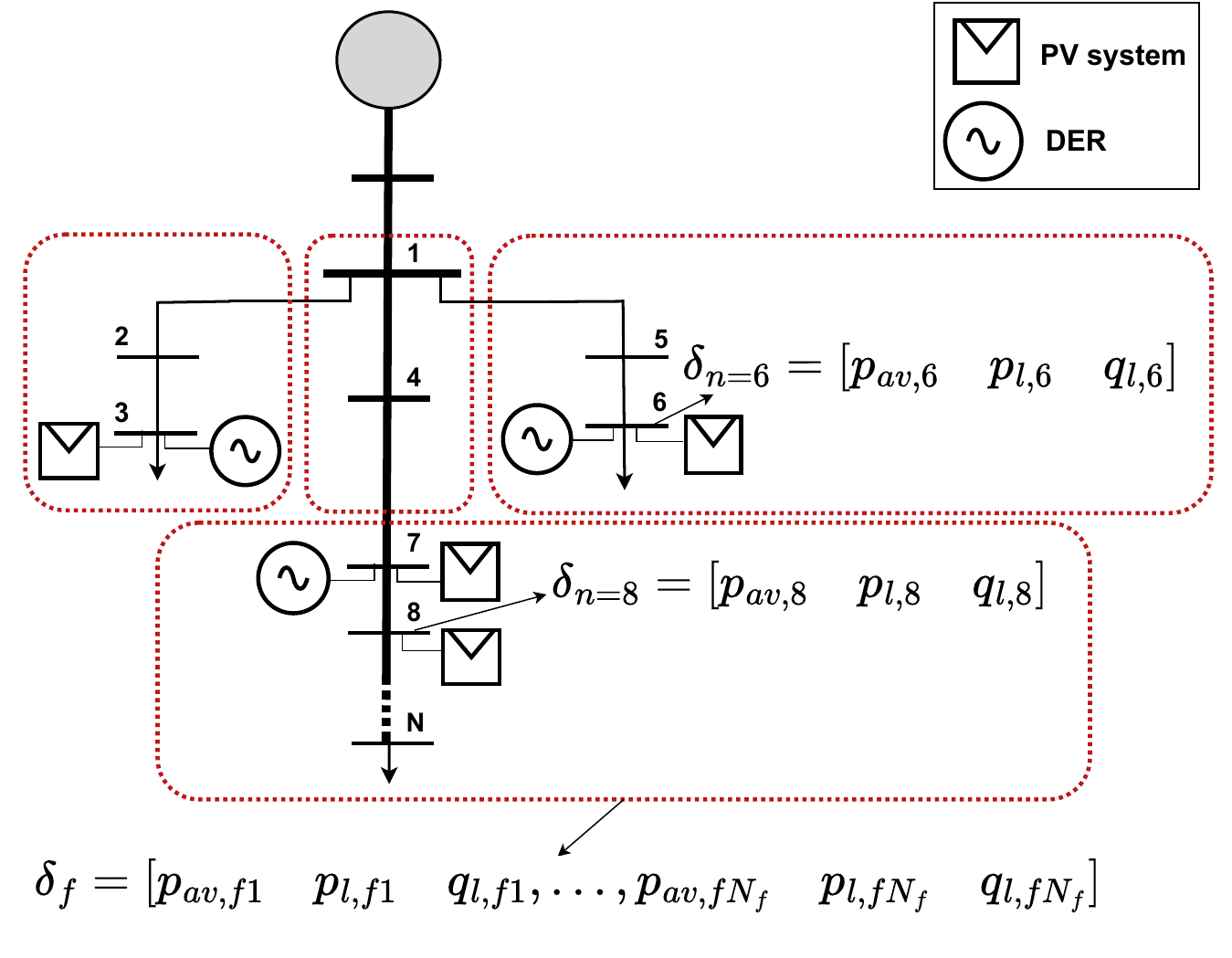}
 \caption{
    Schematic of notations for clusters and nodes. Each node $n$ provides three data points $p_{av,n}$, $p_{l,n}$, $q_{l,n}$. 
    Each data provider $f$  collects data from a ``cluster'' of $N_f$ nodes (indicated by dotted boxes). There are $N$ nodes, $F$ clusters, $3N$ total data points, and $3N_f$ data points in a cluster.}
       \label{fig:schematicofthemodel}
\end{figure}

We assume that the data making up the available samples $\{\widehat{\bm{\delta}}_i\}_{i=1}^I$ is obtained from $F$ individual data providers. 
Each data provider $f=1,...,F$ submits data samples for a subset of $N_f$ nodes.
This modeling choice generalizes the formulation in \cite{mieth2023data}, which modeled a one-to-one relationship between each uncertain parameter and a single data provider. 
This allows us to better reflect real-world organizational structures where multiple data-generating resources or customers are managed by a single service provider, e.g., an aggregator. We call each data stream with data on multiple parameter data \textit{cluster}.
As shown in Fig.~\ref{fig:schematicofthemodel}, we assume that each node exists only in exactly one cluster. 
As a result, $\sum_{f=1}^F N_f = N$ and $\bm{\delta}$ can be rearranged as 
 \begin{align}
\bm{\delta} = \big[
    \bm\delta_1,...,\bm\delta_f,...,\bm\delta_F
    \big]^{\top}
\label{eq:delta_reordered}
\end{align}
where
\begin{flalign}
&\hspace{0em}\bm{\delta}_f\!\!=\!\!\!\big[
    p_{av,f1},\!\!\!\ p_{l,f1},\!\!\
    q_{l,f1}\!,\ \!\!
    ...,\!\!\
    p_{av,f{N_f}},\!\!\!\ p_{l,f{N_f}},\!\!\
    q_{l,f{N_f}}\!\big]^{\!\top}.\label{vectordeltadatacluster1}
\end{flalign}
To avoid notational clutter we will not introduce separate notation for the rearranged $\bm{\delta}$.
Each vector $\bm{\delta}_f\in\mathbb{R}^{3N_f}$ and we assume that each data provider submits $I$ samples for their cluster, i.e.,$\{\widehat{\bm{\delta}}_{f,i}\}_{i=1}^{I}$.  
In the special case of $F=N$, i.e., $N_f = 1,\  f=1,...,F$, the system operator obtains individual datasets from each node.

We note that our model can also accommodate the case where the different parameters of one node are managed by different data providers. 
The special case where each uncertain parameter is informed by an individual data provider recovers the model from \cite{mieth2023data}. 
In this case, we assume that the data making up the available samples $\{\widehat{\bm{\delta}}_i\}_{i=1}^I$ is obtained from $F=3N$ individual data providers. Each data source, $f=1,...,3N$, submits data samples for a specific uncertain parameter related to one node.
As a result, $\bm{\delta}$ can be rearranged as 
 \begin{align}
\bm{\delta} = \big[
    p_{av,1},...,p_{av,N},...,p_{l,1},...,p_{l,N}, q_{l,1},...,q_{l,N}
    \big]^{\top}.
\label{eq:delta_reordered_one-to-one}
\end{align}
Another example is the case of one data provider supplying all data related to loads and one data provider supplying all data on renewable availability. 
In this case $F=2$ and we have $\bm{\delta} = [\bm{\delta_1},\bm{\delta}_2]$ with 
\begin{equation}
\begin{aligned}
    \bm{\delta}_1 &= \{p_{l,1}, q_{l,1},..., p_{l,N}, q_{l,N}\} \\ \bm{\delta}_2 &= \{p_{av,1},..., p_{av,N}\}.
\end{aligned}
\end{equation}
Any organizational structure of the data clusters can be expressed through the arrangement and clustering of $\bm{\delta}$ and only requires the corresponding rearrangement of the rows of $\bm{a}_k$. 

Alongside the dataset $\{\widehat{\bm{\delta}}_{f,i}\}_{i=1}^I$, each data provider submits data quality information $\epsilon_f$, indicating that the unknown true distribution ${\mathbb{P}}_f$ of $\bm{\delta}_f$ is within a Wasserstein distance of at most $\epsilon_f$ from the empirical distribution $\widehat{\mathbb{P}}_f$ supported by $\{\bm{\delta}_{f,i}\}_{i=1}^I$ with high probability.
Intuitively, the DSO can internalize information on distributional ambiguity of the available data into problem \cref{prob:main_op} by restating it as a distributionally robust optimization (DRO) problem:
\begin{subequations}
\begin{align}
& \min_{{\bm{\alpha}, \bm{q}_c,\bm{p}_B,\bm{q}_B,\bm{\delta}}}\quad \sup_{\mathbbmss{Q}\in \mathcal{A}} \mathbb{E}_\mathbbmss{Q}\left[ C\left(\bm{\alpha}, \bm{q_c},\bm{p_B},\bm{q}_B,\bm{\delta}\right)\right] 
\label{eq:msw_dro_objective}
\\
\text{s.t.}\ &\inf_{ \mathbb{Q}\in \mathcal{A}}  \!\!\mathrm{Pr}_{\mathbb{Q}}
\left\{\begin{aligned}
 \max _{k = 1,...,K}[\langle \bm{a}_k, \bm{\delta}\rangle \! +\! c_k ]\le\! \!0 
\end{aligned}\right\} \!\ge \! 1\!\!-\eta^{vol}\label{jointsto}\\
&\inf_{\mathbbmss{Q}\in \mathcal{A}}\!\! \mathrm{Pr}_{\mathbb{Q}}
\left\{\begin{aligned}
\! ((1\!-\!\alpha_n) p_{av,n})^2\!\!+ q_{c,n}^2-S_n^{2} ]\le\! 0 
\end{aligned}\right\}\!\! \ge\! 1-\!\eta^{inv} \nonumber\\& n = 1,...,N \label{eq:thirdddist} \\
&\eqref{eq:deterministiccon3}-\eqref{eq:deterministiccon5}
\end{align}%
\label{prob:msw_dro_problem}%
\end{subequations}%
where similar to $\bm{\delta}$ in \cref{eq:delta_reordered}, we overload the notation of $\bm{a}_k$ to match the order of the entries of $\bm{\delta}$.
Here, unlike in \eqref{prob:main_op}, the distribution $\mathbb{Q}$ is selected from a group of potential distributions known as the ambiguity set $\mathcal{A}$, which is informed by datasets $\{\widehat{\bm{\delta}}_{f,i}\}_{i=1}^{I}$ and their quality $\epsilon_f$ for all $f=1,...,F$.  

To reflect the fact that individual data quality information is available for each dataset, we leverage results from \cite{mieth2023data} and define the multi-source Wasserstein
ambiguity set $\mathcal{A}^{MSW}$
as:
\begin{equation}
\begin{aligned}
    \mathcal{A}^{MSW}=\left\{ \!\mathbb{Q} \!\in\! \mathcal{P}(\Xi) \middle|\!
\begin{array}{l}
\!P_{f\#} \mathbb{Q} = \mathbb{Q}_f,\!\quad  f= 1,..., F,  \\
\!W_p^p(\hat{\mathbb{P}}_f, \mathbb{Q}_f\!) \!\leq \!\epsilon_f, \!\quad   f = 1,..., F,\!\!
\end{array}
\right\}\!,
\end{aligned}
\label{eq:msw_ambiquity_set}
\end{equation}
where $W_p^p(\hat{\mathbbmss{P}}_f, \mathbb{Q}_f)$ denotes the p-Wasserstein distance between the empirical distribution $\hat{\mathbbmss{P}}_f$ supported by $\{\widehat{\bm{\delta}}_{f,i}\}_{i=1}^I$ and $\mathbb{Q}_f$.
Set $\mathcal{P}(\Xi)$ collects all probability distributions defined on support 
$\Xi$. 
The support $\Xi$ may be deduced from the available datasets that are accessible, or it can be established based on various technical factors. 
Moreover, $P_{f\#}$ denotes the push-forward distribution of the joint measure $\mathbb{Q}$ under the projection onto the $f$-th coordinate, i.e.,
\begin{flalign}
    P_{f\#} \mathbb{Q} \coloneqq \!\!\int_{\Xi_{-f}}\hspace{-1.2em}\mathbb{Q}(d\bm\delta_1, \ldots, d\bm\delta_{f-1}, d\bm\delta_f^*, d\bm\delta_{f+1}, \ldots, d\bm\delta_F)
\end{flalign}
with $\Xi_{-f} := \Xi_1 \times ... \times \Xi_{f-1} \times \Xi_{f+1} \times ... \times \Xi_{F} $ and $\Xi_f$ being the projection of $\Xi$ on the $f$-th coordinate. 
As noted in \cite{mieth2023data}, ambiguity set  $\mathcal{A}^{MSW}$
in \cref{eq:msw_ambiquity_set} generalizes established Wasserstein ambiguity sets, e.g., as in \cite{mohajerin2018data}, to accommodate individual Wasserstein budgets $\epsilon_f$ for each component $\bm{\delta}_f$ of $\bm{\delta}$.
The term ``multi-source'' DRO highlights the fact the worst-case distributional shifts of the empirical distributions of subsets of samples are contained a space defined by different Wasserstein distances for each data sources (i.e., data provider).\footnote{The work in \cite{rychener2024wasserstein} uses similar terminology for the study of the distributional risk defined by the intersection of Wasserstein balls.}

The following section derives a data-driven solution approach for \cref{prob:msw_dro_problem} extending upon the established sample-average approach outlined in  Section~II.\ref{sec:data-driven_solution}.

\subsection{Solution with clustered data sources}
Our solution approach builds on the results presented in \cite{mieth2023data}.
Assume a generic DRO problem computing optimal $\bm{x}$ under uncertainty $\bm{\delta}$ for a cost function $c(\bm{x},\bm{\delta})$:
\begin{equation}
    \inf_{\bm{x}\in\set{X}} \sup_{\mathbb{Q}\in\set{A}}\  \mathbb{E}_{\mathbb{Q}}\big[c(\bm{x},\bm{\delta})\big].
\label{eq:generic_DRO}
\end{equation}
With ambiguity set $\set{A} = \set{A}^{MSW}$
as defined in \cref{eq:msw_ambiquity_set}, the generic DRO problem in \cref{eq:generic_DRO} allows for a tractable formulation under the following assumption:
\begin{assumption}\label{as:index_standardization}
All datasets $\{\widehat{\bm{\delta}}_{f,i}\}_{i=1}^I,\ f=1,...,F$ are of identical length $I$ and each index $i=1,...,I$ is standardized so that all samples $\{\widehat{\bm{\delta}}_{f,i}\}_{f=1}^F$ with the same index $i$ can be related, e.g., in terms of a shared time-stamp. (See \cite[Condition C-Std.]{mieth2023data}) 
\end{assumption}
If Assumption~\ref{as:index_standardization} holds, \cref{eq:generic_DRO} can be reformulated as:
\allowdisplaybreaks
\begin{subequations}
\begin{flalign}
&\hspace{0em}\inf_{\bm{x}\in\set{X}, \lambda_f\ge0}\ 
     \sum_{f=1}^F \lambda_f \epsilon_f + \frac{1}{I}\sum_{i=1}^{I} s_{i} \\
&\hspace{0em}\text{s.t.} \quad 
    \!\!\!\!s_{i} \ge \sup_{\bm{\delta}\in\Xi}\ c(\bm{x}, \bm{\delta}) \!- \!\!\sum_{f=1}^F\lambda_f\norm{\bm{\delta}_f -\!\! \widehat{\bm{\delta}}_{f,i}}_p \!\!, i=1,...,I.
\end{flalign}%
\label{eq:mwsdro_general_formulation_stand_data}%
\end{subequations}%
\allowdisplaybreaks[0]%
We refer to \cite[Propositon 2]{mieth2023data} for the detailed proof.
For the DSO problem outlined above, Assumption~\ref{as:index_standardization} is reasonable and we use \cref{eq:mwsdro_general_formulation_stand_data} to bring objective \cref{eq:msw_dro_objective} and chance constraints \cref{jointsto,eq:thirdddist} into a tractable form.

\subsubsection{Chance constraint reformulation}
As in the direct data-driven approach from Section~II.\ref{sec:data-driven_solution} above, we leverage CVaR to reformulate chance constraints \cref{jointsto,eq:thirdddist}.
Enforcing 
\begin{flalign}
\hspace{0em} \!\!\sup_{\mathbb{Q}\in \mathcal{A}} \!\mathbb{Q}CVaR_{\eta^{vol}}
\!\left\{
 [ \max _{k = 1,...,K}\langle \bm{a}_k, \bm{\delta}\rangle \!+\! c_k]\le 0 \!
\right\}\! \ge  \!\!1\!\!-\!\eta^{vol}  \label{eq:voltagecvar1}
\end{flalign}
and
\begin{equation}
\begin{aligned}
\sup_{\mathbb{Q}\in \mathcal{A}} \mathbb{Q}CVaR_{\eta^{inv}}\!\!
\left\{
(1\!-\!\alpha_n)^2 p_{av,n}^2\!+\! q_{c,n}^2\!-\!\!S_n^{2}\le \!0 
\right\}\ge\! 1\!\!-\!\eta^{inv}  \label{eq:thirddcvar}
\end{aligned}
\end{equation}
ensures that \cref{jointsto} and \cref{eq:thirdddist} hold, respectively. 

The resulting reformulations of \cref{eq:voltagecvar1,eq:thirddcvar} are:
\allowdisplaybreaks
\begin{subequations}
\label{eq:maincvar0}
\begin{flalign}
&\sum_{f=1}^{F} \lambda_{f}^{vol} \epsilon_{f} + \frac{1}{I}\sum_{i=1}^{I} s_{i}^{vol} \le \eta^{vol} \varpi^{vol}\label{CVaRvoltagejoint1} \\
& s_{i}^{vol}\! \ge\! c'_k +\!\! \sum_{f=1}^{F}  \sum_{m=1}^{3N_{f}} (z_{k,f} \hat{\delta}_{f,i,m} + u_{k,f,i,m}\overline{\delta}_{f,m} - l_{k,f,i,m} \underline{\delta}_{f,m})\nonumber \\&    \quad k=1,\ldots,K+1 \quad i=1,\ldots,I \\
&a'_{k,f,m} - z_{k,f,i} = u_{k,f,i,m} - l_{k,f,i,m} \quad\forall k,f,i,m\\
& |z_{k,f,i}| \le \lambda_{f}^{vol} \quad\forall k,f,i \\
& u_{k,f,i,m}, l_{k,f,i,m} \ge 0 \quad\forall k,f,i,m \label{CVaRvoltagejoint2} &&
\end{flalign}
\label{eq:CVaRvoltagejoint_reform}
\end{subequations}
\allowdisplaybreaks[0]%
and
\allowdisplaybreaks
\begin{subequations}
\label{eq:maincvarinv}
\begin{flalign}
& \lambda_{n}^{inv} \epsilon_{f(n)} + \frac{1}{I}\sum_{i=1}^{I} s_{n,i}^{inv}  \le  \eta^{inv} \varpi_n^{inv} \quad n = 1,...,N \label{CVaRinvertor1}&& \\
& s_{n,i}^{inv}  \!\ge \!w_{n}^{inv} \!+\! (1\!-\!\!\alpha_{n}\!)^2 [\bm{e}_n^T\overline{\bm{\delta}}_{f(n)}^2]\!  -\! \!\lambda_{n}^{inv} ([\bm{e}_n^T\overline{\bm{\delta}}_{f(n)}] \!- [\!\bm{e}_n^T\!\hat{\bm{\delta}}_{f(n),i}]\!) 
\nonumber\\ & \quad i=1,\ldots,I \quad n = 1,...,N  &&\\
& s_{n,i}^{inv} \ge w_{n}^{inv} \!+ (1-\alpha_{n})^2 [\bm{e}_n^T\hat{\bm{\delta}}_{f(n),i}^2] \quad i=1,\ldots,I\nonumber\\ & n = 1,...,N &&\\
& s_{n,i}^{inv}  \ge 0 \quad i=1,\ldots,I \quad n = 1,...,N.  \label{CVaRinvertor2} && 
\end{flalign}
\label{eq:cvarinverter_reform}%
\end{subequations}%
\allowdisplaybreaks[0]%
We refer to Appendix~\ref{appendix:b} for the detailed derivation steps.
Unlike the reformulation of the voltage chance constraint \cref{eq:CVaRvoltagejoint_reform}, where data from all datasets $\{\bm{\delta}_{f,i}\}_{i=1}^I$ is used, the formulation in \cref{eq:cvarinverter_reform} accounts for the fact that each inverter constraint only requires data related to a single PV system. 
We therefore select a single feature from $\bm{\delta}_f$, expressed as $\bm{e_n}^T\bm{\delta}_f$, where $\bm{e_n}$ is a column vector with all
zero entries except for its n-th index that corresponds to the index of the $p_{av,n}$ in $\bm{\delta}_{f(n)}$ which is equal to one. If a PV is located at node $n$, then $f(n)$ denotes the cluster that includes node $n$.

\subsubsection{Objective function}
We reformulate the objective function using the results of \cite[Proposition~1]{mieth2023data}. To avoid clutter in the main text body we report the final objective in \cref{eq:finalobjective} and additional auxiliary constraints in Appendix~\ref{appendix:c}. 
\subsubsection{Complete formulation}
The final multi-source DRO model formulation is:
\begin{tcolorbox}[mygraybox]
\begin{subequations}
\label{eq:mainfinal55}
\begin{align}
&  \min  \sum_{f=1}^{F} \lambda_f^{co} \epsilon_f +  \sum_{n=1}^{N}\Big( \big(\frac{1}{I} \sum_{i=1}^{I} s_{n,i}^{co1}\big) +  e_n ||q_{c,n}|+|q_{B,n}||\nonumber\\&\qquad  +\big(\frac{1}{I} \sum_{i=1}^{I} s_{n,i}^{co2}\big)  \Big) \label{eq:finalobjective}\\
&\text{s.t.}
\quad\eqref{objfinal1}-\eqref{objfinal2} \quad [\text{Objective aux. constraints}]\nonumber  \\
&\qquad \eqref{eq:deterministiccon3}-\eqref{eq:deterministiccon1}\quad [\text{Deterministic constraints}]\nonumber\\
&\qquad \varpi^{vol} +\varphi^{vol} \le 0 \\
&\qquad \eqref{eq:maincvar0}\quad [\text{Voltage CVaR aux. constraints}]\nonumber\\
&\qquad \varpi_n^{inv} +\varphi_n^{inv} \le 0  \qquad  n = 1,...,N \\
&\qquad \eqref{eq:maincvarinv}\quad [\text{Inverter CVaR aux. constraints}]\nonumber\\
&\qquad  \varphi_n^{inv} \le 0 \qquad n= 1,...,N\\
&\qquad  \varphi^{vol} \le 0
\end{align}
\end{subequations}
\end{tcolorbox}

\subsection{Marginal value of data quality}
Variables $\lambda_f^{co}$, $\lambda_f^{vol}$ and $\lambda_n^{inv} ,n= 1,..., N $ can be used to assess how sensitive the optimal solution of \eqref{eq:mainfinal55} is to changes in the datasets. 
Specifically, each $\lambda_f^{co}$, $\lambda_f^{vol}$ and $\lambda_n^{inv} ,n= 1,..., N $, denotes the maximal variation of the worst-case expected cost in the objective and constraints to marginal changes in the datasets. See also the discussion in \cite{gao2023finite}.
We interpret this as a proxy for how much the optimal solution \textit{relies} on the correctness of a given dataset. 
We can combine these values into a marginal value of data quality $\mu_f$ as follows.
Define $\mathcal{L}$ as the Lagrangian of  \eqref{eq:mainfinal55}, and let $\phi^{vol}$ and $\phi_n^{inv}$ denote the dual multipliers for constraints \eqref{CVaRvoltagejoint1} and \eqref{CVaRinvertor1}. 
For a primal-dual optimal solution \eqref{eq:mainfinal55} we can use the Envelope Theorem to assess the sensitivity of the optimal solution to changes in data quality $\epsilon_f$ submitted by a data provider $f$:
\begin{equation}
\mu_f =\frac{\partial \mathcal{L}}{\partial \epsilon_f} = \lambda^{co}_f + \phi^{vol} \lambda^{vol}_f + \sum_{n \in N} \phi_n^{inv} \lambda_n^{inv}. \label{marginalequation}
\end{equation} 
The marginal value of data quality $\mu_f$  as defined in \cref{marginalequation} consists of three components. The first is the immediate effect of data quality on the uncertain part of the objective, denoted by $\lambda_f^{co}$ the second and third are the effects
of data quality on the objective that occurs indirectly via the cost of enforcing chance constraints \eqref{CVaRvoltagejoint1} and \eqref{CVaRinvertor1}, denoted by $\phi^{vol} \lambda^{vol}_f$ and $\sum_{n \in N_{n,f}} \phi_n^{inv} \lambda_n^{inv}$, respectively. 
Value $\mu_f$ captures the marginal value of increasing data quality in the decision and, as such, offers suitable properties for pricing data similar to marginal pricing approaches for electricity. 
An extensive discussion on the design and properties of a pricing mechanism based on this, however, is subject to further research and beyond the scope of this paper. 

\section{Case Study}
\label{sec:casestudy}
We perform numerical experiments using the multi-source data-driven OPF in \cref{eq:mainfinal55} to assess data utility depending on data quality and operational situation.

\subsection{Data and implementation}
\label{ssec:data_and_implementation}
\begin{figure}
    \centering
    \includegraphics[width=0.98\linewidth]{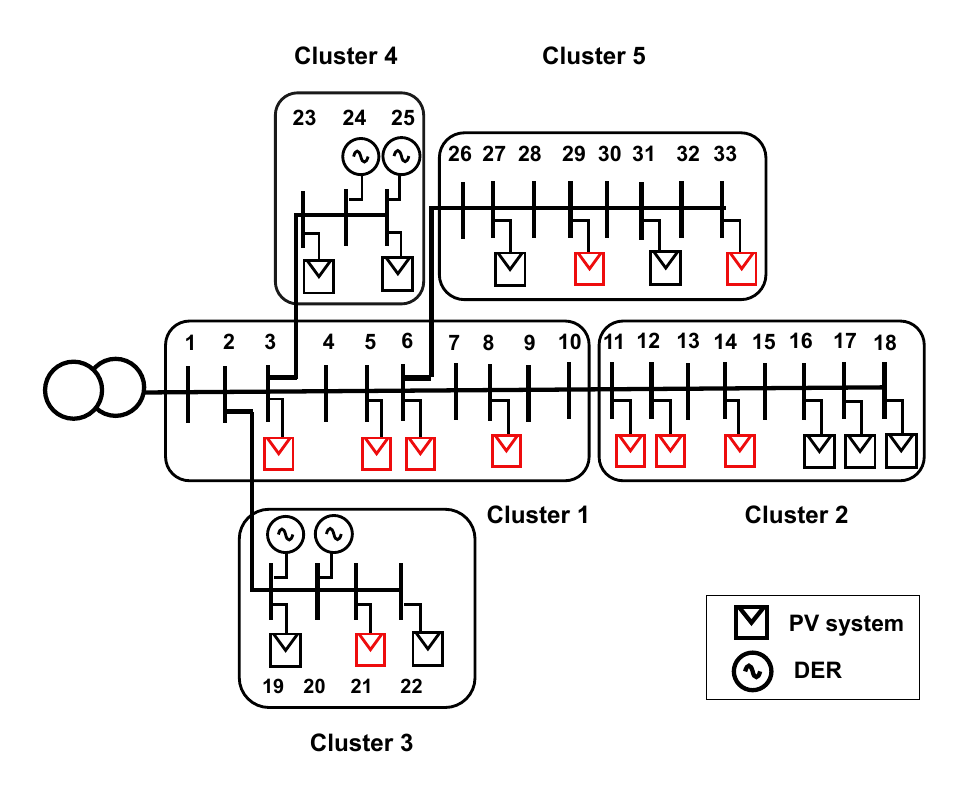}
    \caption{Schematic of the modified IEEE 33-bus system with five data providers, each corresponding to a cluster of nodes. Black PV systems indicate the \textit{Low PV} case; Black + red PV systems indicate the \textit{High PV} case.}
    \label{fig:pv_cases}
\end{figure}

We use the \textit{case33bw} dataset from MATPOWER \cite{matpower2016} as a basis for our test system.
Voltage limits $V_{max}$ and $V_{min}$ are set to 1.1 and 0.9 p.u. (Recall that $V_{max}$, $V_{min}$ limit the square voltage magnitude due to our use of the \textit{LinDistFlow} model. See Section~II.\ref{ssec:problem_formulation} above.)
We modify the original dataset by adding 
9 and 19 PV systems for a \textit{Low PV} and \textit{High PV} scenario. The capacities and placements are indicated in Table~\ref{tab:systemcapacity} and Fig.~\ref{fig:pv_cases}, which also shows system topology. 
We use real PV injection data from the Pecan Street database \cite{PecanStreetDataport} to simulate a 24-hour PV profile as the basis for forecasts $\check{\bm{p}}_{av}$.  
Similarly, we extend the single-period load data from the original \textit{case33bw} dataset to 24 hours by scaling them with load profiles. We consider two load scenarios, \textit{High Load} and \textit{Low Load}. To achieve this, we utilize two load profiles from the ENTSO-E Transparency Platform \cite{OpenPowerSystemDataplatform} and National Grid US \cite{Nationalgridus}, one with high peak load and the other with low peak load (60\% relative to the high peak load profile).
We set $c_n = 10$, $d_n = 3$,
$e_n = 3$ , $h_n = 6$  (following \cite{dall2017chance})
for all $n$ 
and we fix $q_{B,n}$ to zero for all $n$.

For this case study we assume 5 aggregators that act as data providers, i.e., data is coming from 5 clusters with 10, 8, 4, 3, and 8 nodes, respectively, as shown in Fig.~\ref{fig:pv_cases}. 
For every cluster $f = 1,..., 5$ we create a set of $I=25$ samples of $\bm{\delta}_f$ as follows.
For each uncertain parameter $\mathbf{p}_{av}$, $\mathbf{p}_l$, $\mathbf{q}_l$ we use the 24-hour profiles as described above as forecasts $\check{\bm{\delta}} = (\check{\bm{p}}_{av},\check{\mathbf{p}}_l, \check{\mathbf{q}}_l)$.
We then draw the samples $\{\hat{\delta}_{f,i}\}_{i=1}^{I},\ f=1,...,F$ from a normal distribution $N(\check{\bm{\delta}}_f, \Sigma_f)$ with $\Sigma_f = \diag\big((0.2\check{\bm{\delta}}_f)^2\big)$, i.e., the sample distribution is centered around the forecast with a relative standard deviation of 20\% of the forecast value and no correlation.
We highlight that we only use this normal distribution for data generation. It is not 
a prerequisite for the method, which is independent from the underlying forecast error distribution or scenario generation process. 
Finally, the samples are truncated to the support $\Xi_f = [\underline{\bm \delta}_f, \overline{\bm \delta}_f],\ f=1,...,F$.
In practice, support interval $ [\underline{\bm\delta}_f,\overline{\bm\delta}_f]$ can be defined by the DSO, e.g., by inferring it from the available data such that $\underline{\bm\delta}_f = \min_i\{\hat{\bm\delta}_{f,i}\}, \overline{\bm\delta}_f = \max_i\{\hat{\bm\delta}_{f,i}\}$ or in terms of the physical limits of each resource in cluster $f$. In our study we set $\underline{\bm\delta}_f =  \!\!\big[
    0,\!\!\ 0.5\check{p}_{l,f1},\!\!\
    0.5\check{q}_{l,f1},\ \!\!
    ...,\!\!\
    0,\!\ 0.5 \check{p}_{l,f{N_f}},\!\!\
    0.5 \check{q}_{l,f{N_f}}\!\big]^{\!\top}$ and $\overline{\bm\delta}_f =  \!\!\big[
    S_{f1},\!\!\ 1.2\check{p}_{l,f1},\!\!\
    1.2\check{q}_{l,f1},\ \!\!
    ...,\!\!\
    S_{f{N_f}}, 1.2 \check{p}_{l,f{N_f}},\!\!\
    1.2 \check{q}_{l,f{N_f}}\!\big]^{\!\top}.$ 

\begin{remark}
This means, we model the relationship between data and forecasts such that the collection of data from each data provider $\{\hat{\bm\delta}_{i,f}\}_{i=1}^I$ captures the distribution of the forecast errors.
Hence, if the provided forecast $\check{\bm{\delta}}$ has low levels of confidence or is biased, this is reflected in the forecast error distribution.
\end{remark}

All computations have been implemented in the Julia language using JuMP \cite{lubin2023jump} and solved using the SCS solver \cite{ocpb:16}.  They have been performed on a standard MacBook laptop with 16 GB memory and an Apple M3 processor. 
The solving 
time for one time step was around 5 minutes on average. Our implementation is available open source \cite{M.Ghazanfariharandi}.

\begin{table}
\caption{PV capacities and forecasts for 6pm in kW. Nodes with~* have 0 PV capacity in the \textit{Low PV} case.} 
\label{tab:systemcapacity}
\small
\renewcommand{\arraystretch}{1.2}
\setlength{\tabcolsep}{3pt} 
\centering
\begin{tabular}{|c|c|c|c|c|c|c|c|}
\toprule
\textbf{Node} &\textbf{3} & \textbf{5} & \textbf{6} & \textbf{8} &\textbf{ 11} & \textbf{12}& \textbf{14}    \\
\hline
\textbf{Cap.}&500
 &500
 &750
 & 400
&750
&800
&200
\\
\hline
\textbf{$\check{p}_{av,n}$}
& 241.6
& 241.6
& 362.5
& 193.3
& 362.5
&386.6
&96.7
\\
\midrule
\textbf{Node}& \textbf{16\textsuperscript{*}} & \textbf{17\textsuperscript{*}} & \textbf{18\textsuperscript{*}} & \textbf{19\textsuperscript{*}}  & \textbf{21}  & \textbf{22\textsuperscript{*}}& \textbf{23\textsuperscript{*}}\\
\hline
\textbf{Cap.}
&500
&200
&500
&200
&500
&500
&200
\\
\hline
\textbf{$\check{p}_{av,n}$}&241.6
&96.7
&241.6
&96.7
&241.6
&241.6
&96.7
\\
\midrule
\textbf{Node}& \textbf{25\textsuperscript{*}} & \textbf{27\textsuperscript{*}} & \textbf{29} & \textbf{31\textsuperscript{*}}  & \textbf{33}  & \textbf{}&\\
\hline
\textbf{Cap.}
&300
&600
&600
&300
&800
&
&\\
\hline
\textbf{$\check{p}_{av,n}$}&145
&290
&290
&145
&386.6
&
&
\\
\bottomrule
\end{tabular}
\end{table}

\begin{table}[t]
\caption{DER locations and their power limits in kW.}
\label{tab:systemcapacityDER}
\small
\renewcommand{\arraystretch}{1.2}
\setlength{\tabcolsep}{3pt} 
\centering
\begin{tabular}{|c|c|c|c|c|}
\hline
\textbf{Node} &\textbf{19} & \textbf{20} & \textbf{24} & \textbf{25}   \\
\hline
\textbf{Power limit}&50
 &22
 &50
 & 50
\\
\hline
\end{tabular}
\end{table}

\subsection{Decision performance}
\label{ssec:decision_performance}
We first validate the ability of the derived multi-source data-driven OPF in \cref{eq:mainfinal55} to make reliable decisions.
We want to confirm that given a set of samples from each data provider alongside the data quality information the model ensures constraint satisfaction with 
high probability given the unknown, true distribution of the uncertain parameter.
To this end, we select the hour with the highest PV forecast (1pm) and synthesize a large number of $I^{\rm full}=1000$ samples for the uncertain parameters covered by each data cluster using the method described in Section~\ref{ssec:data_and_implementation} above. 
We define the empirical distribution of these $I^{\rm full}$ samples as the true distribution of $\bm{\delta}_f$.
This will enable us to compute the true data quality information for each data provider in the following step, for which
we randomly select $I$ samples from the generated $I^{\rm full}$ samples for each data provider. We then compute the true Wasserstein distance $\epsilon_f$ between the true distribution supported by the $I^{\rm full}$ samples and the empirical distribution supported by the $I$ samples. 
For our specific experiment,
the resulting values for $\epsilon_f$ in \textit{High Load} case 
were $[0.14,0.078,0.021,0.024,0.012]$ and in \textit{Low Load} case were $[0.12, 0.06, 0.022,0.018,0.14]$
. (Note that these values are random as they depend on the randomly selected $I$ samples for each data provider.)
We then solve \cref{eq:mainfinal55} using the selected samples and their true data quality encoded in $\epsilon_f$ to obtain decisions $\bm{\alpha}^*$, $\bm{q}_c^*$, $\bm{p}_B^*$, $\bm{q}_B^*$.
To test the decision, we then draw new $I^{\rm test}=100$ samples from the true distribution, i.e., the set of $I^{\rm full}$ original samples. 
For comparison, we also run \eqref{eq:mainfinal55} modified such that it recovers the sample average approximation (SAA) formulation as in \cite{dall2017chance}, i.e., ignoring any distributional ambiguity formulation. Because the SAA approach effectively assumes that the empirical distribution of the available data samples perfectly captures the underlying distribution of the uncertain parameter, our model in \eqref{eq:mainfinal55} is equivalent to that solution by setting $\epsilon_f=0,\ \forall f$.

\subsubsection{Voltage constraint satisfaction}

\begin{figure*}
\centering
\begin{subfigure}[b]{0.9\textwidth}
\caption{\small \textit{High PV} and \textit{High Load}}
    \label{fig:boxplotvoltagehighload}
    \includegraphics[width=1\linewidth]{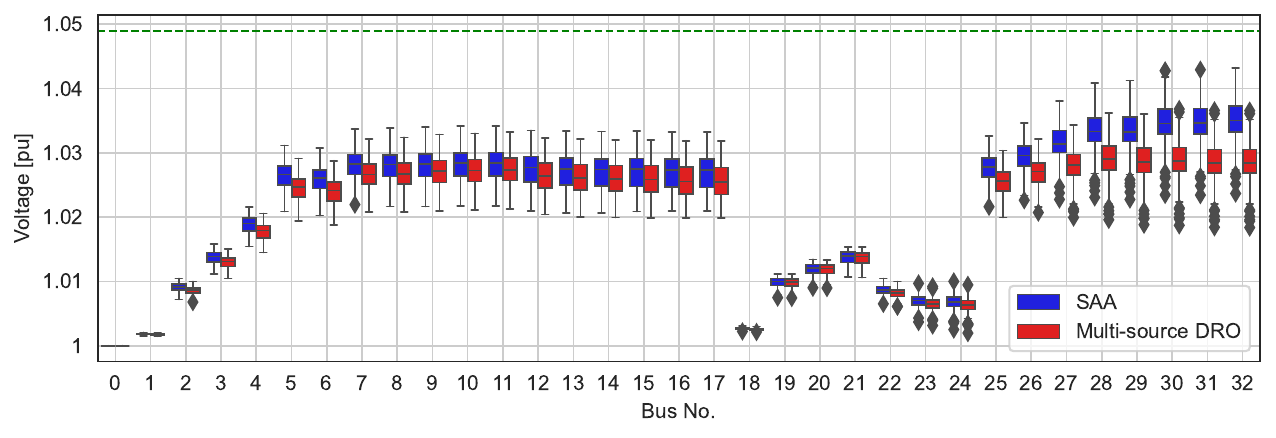}
\end{subfigure}
\begin{subfigure}[b]{0.9\textwidth}
\caption{\small \textit{High PV} and \textit{Low Load}}
    \label{fig:boxplotvoltagelowload}
    \includegraphics[width=1\linewidth]{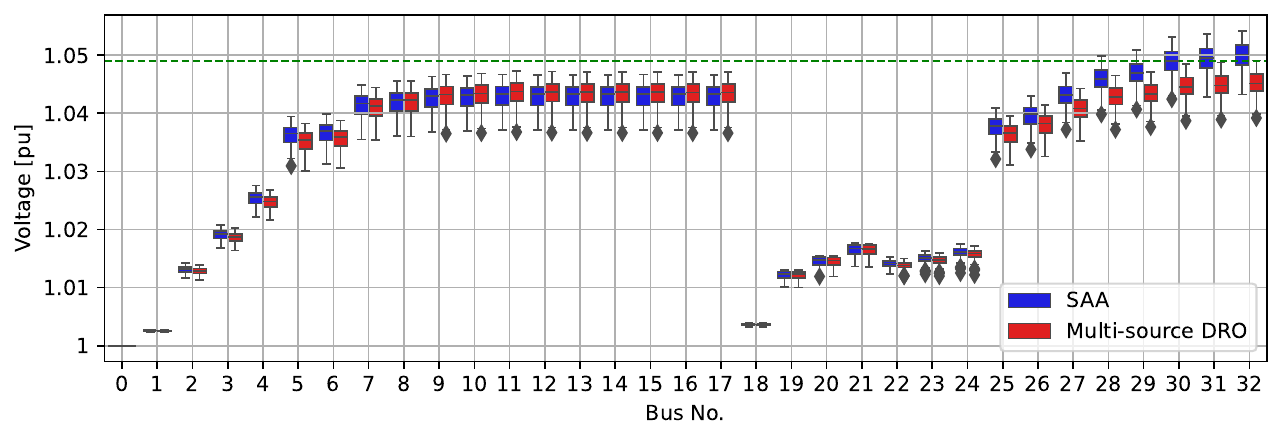}
\end{subfigure}
\caption{
Distribution of voltage magnitudes for the $I^{\rm test}=100$ test scenarios in the \textit{High PV} and \textit{High Load} (a) and \textit{Low Load} (b) cases shown as box plots for the 6pm load. ``SAA'' (blue) shows the out-of-sample performance of the decision if information on data quality is ignored. ``Multi-source DRO'' (red) shows the performance with the correct $\epsilon_f$ for each cluster. The dashed line indicates the upper voltage limit.
For each box plot, the box extends from the 25th to 75th percentile. The center line indicates the median. 
The whiskers span data within 1.5x of the width of the box (inter-quartile range) from the box. Diamonds denote outliers.
}
\label{fig:boxplotvoltage}
\end{figure*}

Fig.~\ref{fig:boxplotvoltage} shows the resulting voltages for the out-of-sample scenarios for the two load cases (\textit{High Load} in Fig.~\ref{fig:boxplotvoltagehighload} and \textit{Low Load} in Fig.~\ref{fig:boxplotvoltagelowload}). 
In the \textit{High Load} case, we observe that both the SAA approach and the proposed multi-source DRO approach ensure safe voltages limits. 
In the \textit{Low Load} case, however, the SAA approach fails to ensure the target risk level of $\eta^{\rm vol} = 5\%$ and at least one voltage limit is violated in $67\%$ of the cases. 
For the multi-source DRO approach, the empirical voltage violation risk is $1\%$. 
(In the \textit{Low Load} case, voltages are more susceptible to PV uncertainty due to the lower local consumption at the grid edge.)
To further confirm the results for the \textit{High Load} case, we also repeated our out-of-sample experiment with higher levels of uncertainty. As we can observe in Table~\ref{tab:outofsamplehighload}, 
in this case some voltage constraint violations  were 
observed for the test cases, but the empirical probability of constraint violation remained within the desired risk level $\eta^{\rm vol} = 5\%$.

\subsubsection{System cost}
\begin{figure}
    \centering
    \includegraphics[width=0.9\linewidth]{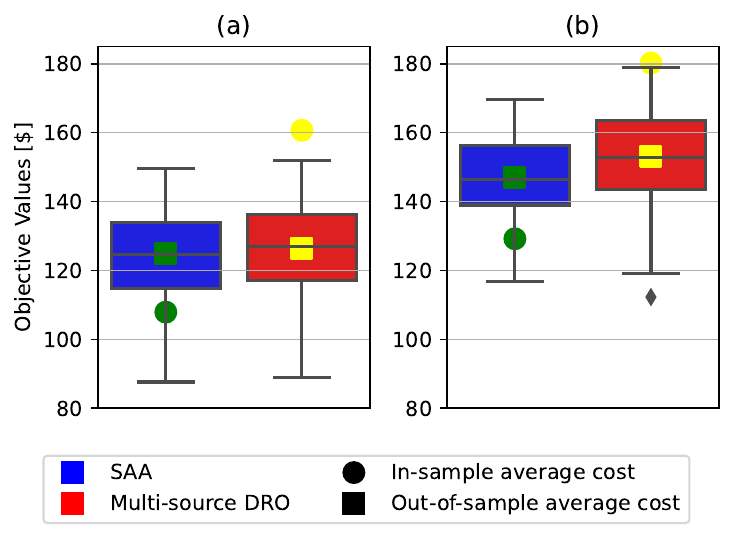}
    \caption{
    Out-of-sample cost distribution for SAA (blue) and multi-source DRO (red) approaches shown as box plots for the \textit{High PV} and \textit{Low Load} (a) and \textit{High Load} (b) cases at 6pm. Circles and squares indicate the corresponding in-sample and out-of-sample cost.
    See caption of Fig.~\ref{fig:boxplotvoltage} for more information on box plots.}
    \label{fig:boxplotobj}
\end{figure}

Fig.~\ref{fig:boxplotobj} shows the resulting cost for the studied cases. 
For both the \textit{Low Load} (Fig.~\ref{fig:boxplotobj}a) and the \textit{High Load} (Fig.~\ref{fig:boxplotobj}b) cases the average out-of sample cost are lower than the worst-case expected cost computed by the multi-source DRO approach.
Hence, the proposed approach reliably protects the DSO from underestimating the expected cost. In other words, the DRO approach successfully computes an upper bound on the expected cost. See also \cite{mohajerin2018data} for a discussion on this property of DRO.
The SAA solution, on the other hand, consistently leads to higher cost out-of-sample, thus systematically underestimating the risk for the decision-maker.
The out-of-sample cost of the multi-source DRO solution in the \textit{Low Load} case is only 1.3\% higher than the SAA but effectively reduces the risk of voltage limit violations as discussed above.
Similarly, for the \textit{High Load} case, the DRO solution provides improved risk quantification and additional robustness at a 4.7\% cost increase.

\begin{table}
\caption{Empirical out-of-sample violation probability of the voltage joint chance constraint in the \textit{High Load} case. }
\label{tab:outofsamplehighload}
\small
\renewcommand{\arraystretch}{1.2}
\setlength{\tabcolsep}{3pt} 
\centering
\begin{tabular}{|>{\columncolor{lightgray!30}}c!{\vrule width 1.5pt}c|c|c|c|}
\rowcolor{lightgray!30} 
\hline
\textbf{Relative standard deviation} 
&\textbf{0.2} 
& \textbf{0.5} 
& \textbf{0.85} 
& \textbf{0.9}   \\
\noalign{\hrule height 1.5pt}
\hline
\textbf{Violation probability}&0\%
 &0\%
 & 2\%
 & 3\%
\\
\hline
\end{tabular}
\end{table}

\subsection{Data value analysis.}

This section demonstrates how our method can be used to quantify data value. For each pair of
PV scenario (\textit{Low PV} and \textit{High PV}) and load scenario (\textit{Low Load} and \textit{High Load}) 
we solve model \cref{eq:mainfinal55} 24 times (once for each hour of the day) using the respective load and PV forecasts generated as described above.
In contrast to the experiment in Section~\ref{ssec:decision_performance} above, we now want to study how data quality impacts the decision and the data value. 
To this end, we conduct this experiment multiple times with predefined data quality values, i.e., 
1.0, 0.1,0.01, 0.005, 0.001, and 0.0001.
For now, we assume that each provider submits data of the same quality.

\subsubsection{Voltage and Objective}
\begin{figure}
    \centering
    \includegraphics[width=1\linewidth]{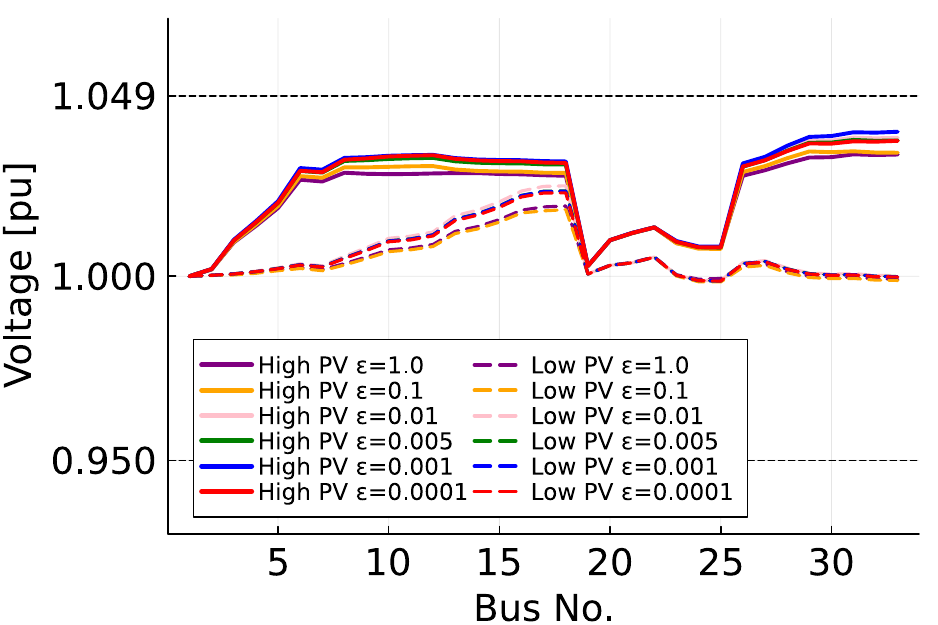}
  \caption{Voltage magnitude profiles across all nodes under varying data quality conditions, showcasing two distinct cases in the \textit{High Load} scenario
  :  \textit{Low PV} case and \textit{High PV} case. Horizontal dashed lines indicate upper and lower voltage limits.}
     \label{fig:voltagecomparisonfig}
\end{figure}

\begin{figure}
    \centering
    \includegraphics[width=1\linewidth]{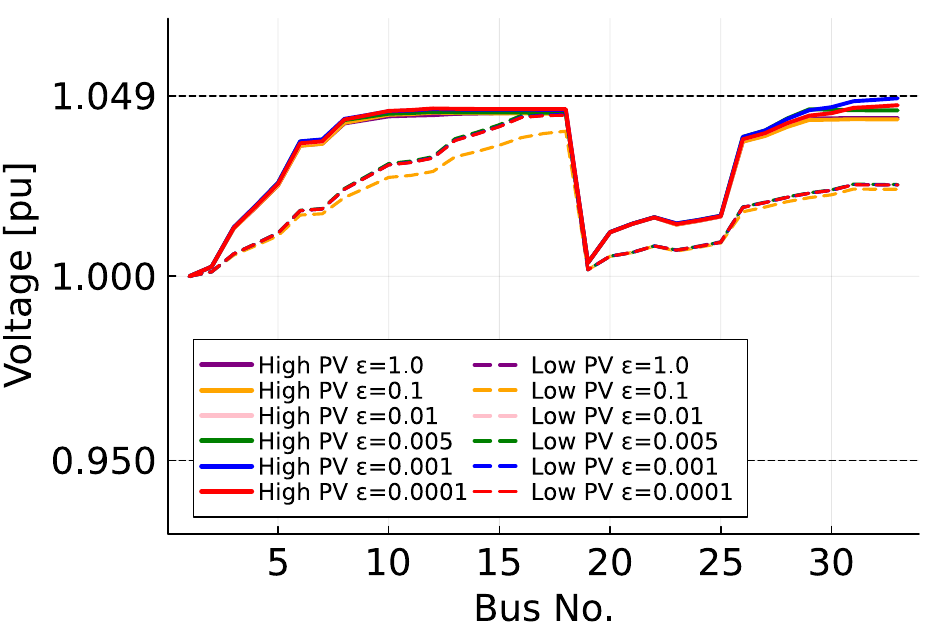}
    \caption{ Voltage magnitude profiles across all nodes under varying
data quality conditions, showcasing two distinct cases in the \textit{Low Load} scenario: \textit{Low PV} case
and \textit{High PV} case. Horizontal dashed lines indicate upper and lower
voltage limits.}
    \label{fig:voltagecomparison_R1}
\end{figure}

Fig.~\ref{fig:voltagecomparisonfig} and Fig.~\ref{fig:voltagecomparison_R1}
show the resulting (expected) voltage profile across the nodes and Table~\ref{tab:objectivevalue} and Table~\ref{tab:objectivevalue_R1} 
show the resulting objective value.
These two tables
show the
worst-case expected cost decreases
with improved data quality. The higher renewable curtailments in the \textit{High PV} case result in higher objective values compared to the \textit{Low PV} case. 
Also, because of the lower total power demand, the cost of the system in the \textit{Low Load} case is systematically lower than in the \textit{High Load} case.
We also observe in Fig.~\ref{fig:voltagecomparisonfig} and Fig.~\ref{fig:voltagecomparison_R1}
that with lower data quality the model chooses a lower voltage baseline to accommodate more PV in the \textit{High PV} scenario. In addition, we can observe that in the \textit{Low Load} case the voltage level is generally higher than in the \textit{High Load} case due to the increased influence of PV injections on the voltages.

\begin{table}
 \caption{Objective values [\$] for various $\epsilon_f$ at 6pm in \textit{High Load} case }
\label{tab:objectivevalue}
\small
\renewcommand{\arraystretch}{1.2}
\setlength{\tabcolsep}{4pt} 
\centering
\begin{tabular}{|>{\columncolor{lightgray!30}}c!{\vrule width 1.5pt}c|c|c|c|c|c|}
\hline
\multirow{2}{*}{} & \multicolumn{6}{c|}{\cellcolor{lightgray!30}$\bm \epsilon_f$} \\
\cline{2-7}
\rowcolor{lightgray!30} 
&\textbf{1} & \textbf{0.1} & \textbf{0.01} & \textbf{0.005} &\textbf{ 0.001} & \textbf{0.0001} \\
\noalign{\hrule height 1.5pt}
\hline
\textbf{High PV}&285.5
 &209.6
 &171.5
 & 169
&166.6
&166.1
\\
\hline
\textbf{Low PV}& 204.4
& 186.7
& 148
& 145.4
& 143.4
&142.9
\\
\hline
\end{tabular}
\end{table}

\begin{table}
 \caption{ Objective values [\$] for various $\epsilon_f$ at 6pm in \textit{Low Load} case} 
\label{tab:objectivevalue_R1}
\small
\renewcommand{\arraystretch}{1.2}
\setlength{\tabcolsep}{4pt} 
\centering
\begin{tabular}{|>{\columncolor{lightgray!30}}c!{\vrule width 1.5pt}c|c|c|c|c|c|}
\hline
\multirow{2}{*}{} & \multicolumn{6}{c|}{\cellcolor{lightgray!30}$\bm \epsilon_f$} \\
\cline{2-7}
\rowcolor{lightgray!30}
&\textbf{1} & \textbf{0.1} & \textbf{0.01} & \textbf{0.005} &\textbf{ 0.001} & \textbf{0.0001} \\
\noalign{\hrule height 1.5pt}
\hline
\textbf{High PV}&192.09
 &132.86
 &110.5
 &108.8
&106.05
&105.2
\\
\hline
\textbf{Low PV}& 144.6
& 100.9
& 84.9
& 82.5
& 80.7
& 79.6
\\
\hline
\end{tabular}
\end{table}

\subsubsection{Data value per constraint}
Tables~\ref{tab:lambdafirstandsecondhigh} and ~\ref{tab:lambdathirdhigh} show the resulting values of $\lambda_{f}^{co}$, $\lambda_{f}^{vol}$, and \mbox{$\lambda_n^{inv}$ , $n = 1,...,N$}, respectively, exemplary for the \textit{High PV} and \textit{High Load }
scenario at 6pm.
We chose this hour because it generally produces the highest $\mu_f$ across all $f=1,...,F$, highlighting the impact of changes in data quality on the decision most prominently.
These values allow us to study data usefulness for each data provider $f$ \textit{and} the various constraints.
When $\epsilon_f$ is set to $1$ for all $f$, i.e., data quality is low, we observe that $\lambda_{f}^{vol} = 0$ and $\lambda_{f}^{co} = 0$ except in clusters 1 and 2, that and $\lambda_n^{inv}$ $\forall n=1,...,N$ are small. This indicates inadequate data quality for clusters 3, 4 and 5, i.e., data quality is so low that the solution does not rely on the data at all. 
The small but non-zero values of $\lambda^{co}$ and $\lambda^{inv}$ for clusters 1 and 2, on the other hand, indicate that even low-quality data will still support the decision.
Decreasing $\epsilon_f$, i.e., improving data quality, implies smaller required ambiguity sets around the empirical data distributions, allowing the operator to place greater trust in the data and reduce uncertainty surrounding forecast values. This leads to reduced system costs (as observed in Table~\ref{tab:objectivevalue}), and an increase in  $\lambda_{f}^{vol}$, $\lambda_{n}^{inv} , n=1,...,N$, and $\lambda_f^{co}$. 
The results in Tables~\ref{tab:lambdafirstandsecondhigh} and ~\ref{tab:lambdathirdhigh} also show that the relationship between data quality $\epsilon_f$ and $\lambda_{f}^{vol}$, $\lambda_f^{co}$, and $\lambda_{n}^{inv}$ does not strictly adhere to a monotonic progression, e.g., in the transition
from $\epsilon_f = 0.001$ to $\epsilon_f = 0.0001$, $ f = 1,...,5$. 

\begin{table}
 \caption{Results for $\lambda_f^{vol}$, $\lambda_f^{co}$
for various $\epsilon_f$ in the \textit{High PV} and \textit{High Load} case.}
\label{tab:lambdafirstandsecondhigh}
\small
\renewcommand{\arraystretch}{1.2}
\setlength{\tabcolsep}{4pt}
\centering
\begin{tabular}{|>{\columncolor{lightgray!30}}c!{\vrule width 1.5pt}c|c|c|c|c|c!{\vrule width 1.5pt}>{\columncolor{lightgray!30}}c|}
\hline
\multirow{1}{*}{\textbf{f}} & \multicolumn{6}{c!{\vrule width 1.5pt}}{\cellcolor{lightgray!30}$\bm \epsilon_f$} & \\
\cline{2-7}
\tabucline{2-7}
\rowcolor{lightgray!30} 
& \textbf{1} & \textbf{0.1} & \textbf{0.01} & \textbf{0.005} &\textbf{ 0.001} & \textbf{0.0001} & \multirow{2}{*}{} \\
\noalign{\hrule height 1.5pt}
\textbf{1} & 0 & 0 & 0.0073 & 0.0073 & 0.0074 & 0.0073 & \\
\textbf{2} & 0 & 0.0052 & 0.04 & 0.04 & 0.04 & 0.04 & \\
\textbf{3} & 0 & 0 & 0.0011 & 0.0012 & 0.0012 & 0.0012 & \\
\textbf{4} & 0 & 0 & 0 & 0 & 0 & 0 & \\
\textbf{5} & 0 & 0 & 0.017 & 0.03 & 0.04 & 0.044 & \multirow{-5}{*}{$\bm \lambda_f^{vol}$} \\
\noalign{\hrule height 1.5pt}
\textbf{1} & 1.34 & 4.06 &9.96 & 9.93 & 10.004 & 19.48 & \\
\textbf{2} & 1.8 & 4.78 & 9.99 & 9.99 & 10 & 15.28 & \\
\textbf{3} & 0 &3& 10.03 & 10.01 & 10& 10.1 & \\
\textbf{4} & 0 & 10 & 10 & 10 & 10 &19.86 & \\
\textbf{5} & 0 & 3.8 & 9.97 & 9.98 & 10.004 & 18.42 & \multirow{-5}{*} {$\bm \lambda_f^{co}$} \\
\hline
\end{tabular}
\end{table}

\begin{table}
 \caption{ Results for $\lambda_n^{inv}$
for various $\epsilon_f$ in the \textit{High PV} and \textit{High Load} case.}
\label{tab:lambdathirdhigh}
\small
\renewcommand{\arraystretch}{1.2}
\setlength{\tabcolsep}{4pt} 
\centering
\begin{tabular}{|>{\columncolor{lightgray!30}}c!{\vrule width 1.5pt}>{\columncolor{lightgray!30}}c!{\vrule width 1.5pt}c|c|c|c|c|c|}
\hline
\multirow{1}{*}{\textbf{f}} & \multirow{1}{*}{\textbf{n}} & \multicolumn{6}{c|}{\cellcolor{lightgray!30}$\bm \epsilon_f$} \\
\cline{3-8}
\rowcolor{lightgray!30} 
& &  \textbf{1} & \textbf{0.1} & \textbf{0.01} & \textbf{0.005} &\textbf{ 0.001} & \textbf{0.0001} \\
\noalign{\hrule height 1.5pt}
\textbf{1}& \textbf{3}&0 &0 &0 & 0&0.034 &0.034 \\
\textbf{1}& \textbf{5}& 0.0001& 0& 0& 0.047& 0.047& 0.047\\
\textbf{1}& \textbf{6}& 0.0003& 0.0006& 0& 0.053& 0.054& 0.054\\
\textbf{1}& \textbf{8}& 0& 0& 0& 0& 0& 0\\
\hline
\textbf{2}&\textbf {11}& 0.0028& 0.002& 0&0.0027& 0.057& 0.057\\
\textbf{2}&\textbf{12}& 0.0004& 0.0014& 0.4&0.083&0.085& 0.084\\
\textbf{2}&\textbf{14} & 0& 0.0001& 0& 0&0.02& 0.019\\
\textbf{2}&\textbf{16} & 0.0004& 0.0008& 0.0002& 0& 0.0003& 0.0002\\
\textbf{2}&\textbf{17} & 0.0001&0.0002& 0& 0& 0.024& 0.021\\
\textbf{2}&\textbf{18} & 0.0018& 0.0006& 0.035& 0.26& 0.036& 0.035\\
\hline
\textbf{3}&\textbf{19} & 0& 0& 0& 0& 0.015& 0.015\\
\textbf{3}&\textbf{21} & 0& 0& 0& 0& 0& 0\\
\textbf{3}&\textbf{22} & 0& 0& 0& 0.075& 0.077& 0.077\\
\hline
\textbf{4}&\textbf{23} & 0& 0& 0& 0& 0.009& 0.009\\
\textbf{4}&\textbf{25} & 0& 0& 0& 0.017& 0.017& 0.017\\
\hline
\textbf{5}&\textbf{27} & 0& 0& 0
& 0& 0& 0\\
\textbf{5}&\textbf{29} & 0& 0& 0& 0.0002& 0& 0\\
\textbf{5}&\textbf{31} & 0& 0& 0& 0& 0& 0\\
\textbf{5}&\textbf{33} & 0.0004& 0.0053& 0& 0.0024& 0.08& 0.081\\
\hline
\end{tabular}
\end{table}
\subsubsection{Data value per data provider}
\begin{figure}
    \centering
    \includegraphics[width=0.9\linewidth,
    trim={0cm 4.1cm 0cm 0cm},clip]{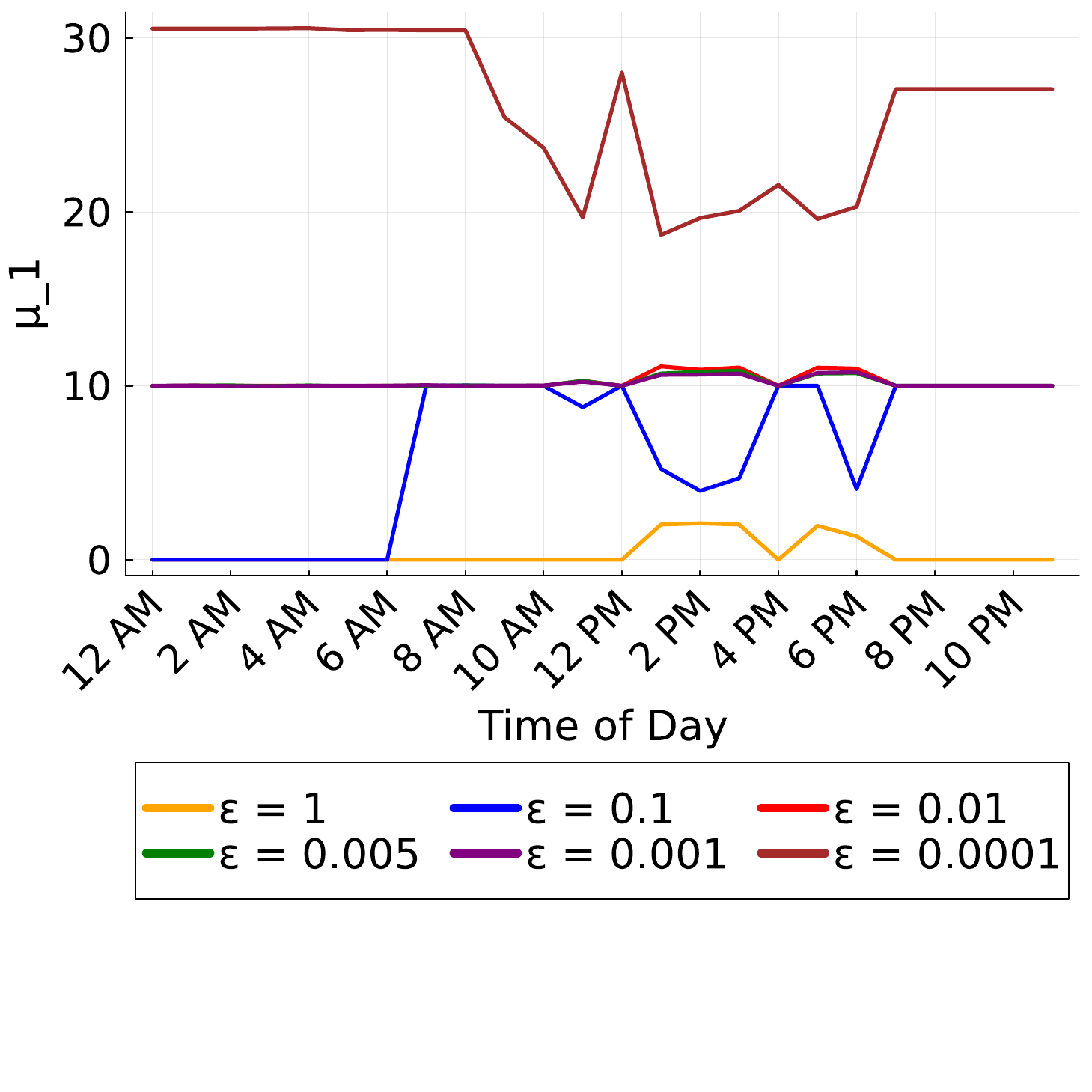}
    \caption{Sensitivity of the optimal objective value to changes of the data quality submitted by cluster 1 over 24 hours of the day in the \textit{High PV} case.}
    \label{fig:mu_1sensitivity}
\end{figure}

\begin{figure}
    \centering
    \includegraphics[width=0.9\linewidth]{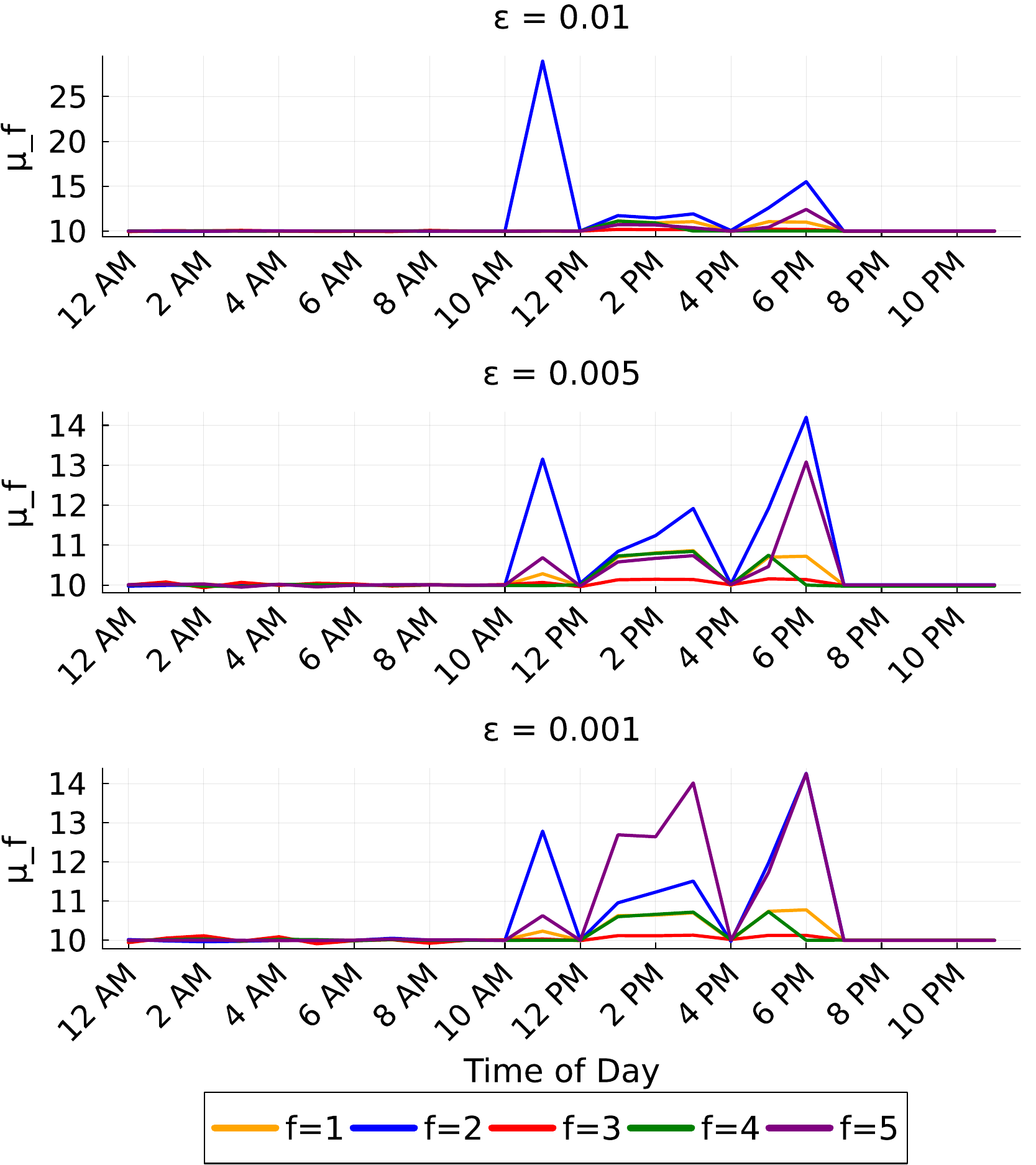}
    \caption{$\mu_f, f= 1,...,5$ over 24 hours of the day in different $\epsilon_f$ in the \textit{High PV} case.}
   \label{fig:threeplot}
\end{figure}

We now turn towards the marginal value of data quality for each data provider, expressed by $\mu_f$ as analyzed in \eqref{marginalequation}.
Tables~\ref{tab:mu_high}, ~\ref{tab:mu_low}, ~\ref{tab:mu_high_lowload_R1}, and ~\ref{tab:mu_low_lowload_R1} show these values for the \textit{High PV} and \textit{Low PV} cases in \textit{High Load} and \textit{Low Load} scenarios
, respectively. 
Figs.~\ref{fig:mu_1sensitivity} and \ref{fig:threeplot} visualize $\mu_f$ over the time of day. 
In Fig.~\ref{fig:threeplot} we observe that the marginal data value of providers 2 and 5 often dominates the others, which we explain with the high total PV capacity in their respective node clusters. 
Fig.~\ref{fig:mu_1sensitivity} provides details on Cluster 1, showing $\mu_1$ for 6 different levels of data quality.
For $\epsilon_1=1$, marginal data quality is consistently low across all times of the day, indicating that the optimal decision relies very little or not at all on the provided data at this quality level.
As $\epsilon_1$ decreases, i.e., data quality increases, the curves tend to have more variation and higher marginal data quality values, suggesting, as expected, that the optimal decision relies more on the information from the provided data if it is \mbox{of higher quality}.

\begin{table}
\caption{Results for $\mu_f$ for various $\epsilon_f$ in the \textit{High PV}  and \textit{High Load} case.}
\label{tab:mu_high}
\small
\renewcommand{\arraystretch}{1.2}
\setlength{\tabcolsep}{4pt}
\centering
\begin{tabular}{|>{\columncolor{lightgray!30}}c!{\vrule width 1.5pt}c|c|c|c|c|c|}
\hline
\multirow{1}{*}{\textbf{f}} & \multicolumn{6}{c|}{\cellcolor{lightgray!30}$\bm \epsilon_f$} \\
\cline{2-7}
\rowcolor{lightgray!30} 
& \textbf{1} & \textbf{0.1} & \textbf{0.01} & \textbf{0.005} & \textbf{0.001} & \textbf{0.0001} \\
\noalign{\hrule height 1.5pt}
\textbf{1} & 1.34 & 4.06  & 10.99  & 10.72  & 10.77 & 20.3 \\
\textbf{2} & 1.8  & 5.45  & 15.5   & 14.2   & 14.26 & 19.78 \\ 
\textbf{3} & 0    & 3     & 10.19  & 10.13  & 10.12 & 10.23 \\ 
\textbf{4} & 0    & 9.99  & 10     & 10     & 10    & 19.86 \\ 
\textbf{5} & 0    & 3.8   & 12.41  & 13.08  & 14.25 & 23.26 \\
\hline
\end{tabular}
\end{table}

\begin{table}[t]
\caption{Results for $\mu_f$ for various $\epsilon_f$ in the \textit{Low PV} and \textit{High Load} case.}
\label{tab:mu_low}
\small
\renewcommand{\arraystretch}{1.2}
\setlength{\tabcolsep}{4pt}
\centering
\begin{tabular}{|>{\columncolor{lightgray!30}}c!{\vrule width 1.5pt}c|c|c|c|c|c|}
\hline
\multirow{1}{*}{\textbf{f}} & \multicolumn{6}{c|}{\cellcolor{lightgray!30}$\bm \epsilon_f$} \\
\cline{2-7}
\rowcolor{lightgray!30} 
& \textbf{1} & \textbf{0.1} & \textbf{0.01} & \textbf{0.005} &\textbf{ 0.001} & \textbf{0.0001} \\
\noalign{\hrule height 1.5pt}
\textbf{1}& 0    & 0     & 9.96   & 10    & 10.02  & 29.64  \\
\textbf{2}& 0    & 2.99  & 10.23  & 10.2  & 10.12  & 23.02  \\ 
\textbf{3}& 0    & 3.08  & 9.97   & 10    & 10     & 16.77  \\ 
\textbf{4}& 0    & 10.01 & 10.01  & 10.01 & 10.01  & 22.68  \\ 
\textbf{5}& 0    & 3.9   & 10.04  & 10.02 & 10.02  & 23.41  \\
\hline
\end{tabular}
\end{table}

\begin{table}
\caption{ Results for $\mu_f$ for various $\epsilon_f$ in the \textit{High PV} and \textit{Low Load} case.}
\label{tab:mu_high_lowload_R1}
\small
\renewcommand{\arraystretch}{1.2}
\setlength{\tabcolsep}{4pt}
\centering
\begin{tabular}{|>{\columncolor{lightgray!30}}c!{\vrule width 1.5pt}c|c|c|c|c|c|}
\hline
\multirow{1}{*}{\textbf{f}} & \multicolumn{6}{c|}{\cellcolor{lightgray!30}$\bm \epsilon_f$} \\
\cline{2-7}
\rowcolor{lightgray!30} 
& \textbf{1} & \textbf{0.1} & \textbf{0.01} & \textbf{0.005} &\textbf{ 0.001} & \textbf{0.0001} \\
\noalign{\hrule height 1.5pt}
\textbf{1}& 2.003 &5.27&4.21  &12.18  & 11.97  & 40.42 \\
\textbf{2}&5.84  &6.33&6.04   &6.043  &14.17   &25.54  \\ 
\textbf{3}&0     &3.29&10.05  &10.09  &10.07   &34.14  \\ 
\textbf{4}&0     &2.67&10.48  &10.62  &10.6    &44.47  \\ 
\textbf{5}&0.067 &5.94&6.01   &12.97  &10.13   &36.15  \\
\hline
\end{tabular}
\end{table}

\begin{table}[t]
\caption{Results for $\mu_f$ for various $\epsilon_f$ in the \textit{Low PV} and \textit{Low Load} case.}
\label{tab:mu_low_lowload_R1}
\small
\renewcommand{\arraystretch}{1.2}
\setlength{\tabcolsep}{4pt}
\centering
\begin{tabular}{|>{\columncolor{lightgray!30}}c!{\vrule width 1.5pt}c|c|c|c|c|c|}
\hline
\multirow{1}{*}{\textbf{f}} & \multicolumn{6}{c|}{\cellcolor{lightgray!30}$\bm \epsilon_f$} \\
\cline{2-7}
\rowcolor{lightgray!30}
& \textbf{1} & \textbf{0.1} & \textbf{0.01} & \textbf{0.005} &\textbf{ 0.001} & \textbf{0.0001} \\
\noalign{\hrule height 1.5pt}
\textbf{1} & 0    & 0    & 4.97   & 10.78  & 9.81   & 36.91  \\
\textbf{2} & 0    & 6.12 & 6.28   & 11.05  & 11.82  & 25.1   \\ 
\textbf{3} & 0    & 3.06 & 10.06  & 10.01  & 10.003 & 18.11  \\ 
\textbf{4} & 0    & 2.3  & 10.48  & 11.84  & 10.4   & 21.62  \\ 
\textbf{5} & 0    & 3    & 3      & 15.09  & 12.73  & 29.94  \\
\hline
\end{tabular}
\end{table}

The marginal value of data quality for each data provider also depends on the data quality of \textit{other} data providers. To highlight this, we conducted an experiment where we set $\epsilon_f = 0.01$ for all $f = 1,3,4,5$ (i.e., all except $2$) varied only $\epsilon_2$ to observe its effect on the marginal values. 
Table~\ref{tab:mu_low_cluster2changes} itemizes these results. 
First, as discussed above, we observe that the marginal value of data quality does not increase monotonically.
Increasing data quality will increase the model's reliance on the provided data, leading to a general increase in
marginal data quality value. On the other hand, improved data will improve the decision, i.e., reducing the
cost of
some marginal cost components that contribute to $\mu_f$ as per \cref{marginalequation}.

\begin{table}
\caption{Results for $\mu_f$ for various $\epsilon_2$ and $\epsilon_f = 0.01$, $f=1,3,4,5$ in the \textit{High PV}  and \textit{High Load} case.}
\label{tab:mu_low_cluster2changes}
\small
\renewcommand{\arraystretch}{1.2}
\setlength{\tabcolsep}{4pt}
\centering
\begin{tabular}{|>{\columncolor{lightgray!30}}c!{\vrule width 1.5pt}c|c|c|c|c|c|}
\hline
\multirow{1}{*}{\textbf{f}} & \multicolumn{6}{c|}{\cellcolor{lightgray!30}$\bm \epsilon_2$} \\
\cline{2-7}
\rowcolor{lightgray!30} 
& \textbf{1} & \textbf{0.1} & \textbf{0.01} & \textbf{0.005} &\textbf{ 0.001} & \textbf{0.0001} \\
\noalign{\hrule height 1.5pt}
\textbf{1} & 10.61 & 11.57 & 10.99  & 10.78 & 10.68  & 10.47  \\
\cellcolor{gray!50}\textbf{2} & \cellcolor{gray!50}1.79 & \cellcolor{gray!50}5.77 & \cellcolor{gray!50}15.5 & \cellcolor{gray!50}14.13 & \cellcolor{gray!50}13.58  & \cellcolor{gray!50}17.95 \\ 
\textbf{3} & 10.1  & 10.17  & 10.19  & 10.17  & 10.11  & 10.08  \\ 
\textbf{4} & 10    & 10     & 10     & 9.99   & 10     & 9.98   \\ 
\textbf{5} & 11.4  & 9.96   & 12.41  & 11.83  & 11.67  & 12.3   \\
\hline
\end{tabular}
\end{table}

\subsubsection{Critical data quality}
The system operator can apply the proposed method to identify the minimal level of data quality required from each data provider so that the data offers some utility to the decision-making process. We define:
\begin{definition}[Critical $\epsilon_f$]
    The smallest value $\epsilon_f$ for which any $\bm{\lambda}^{vol}$, $\bm{\lambda}^{inv}$, $\bm{\lambda}^{co}$ associated with data provider $f$ is equal to zero is called \emph{critical} $\epsilon_f$ with respect to a given constraint.
\end{definition}
We refer also to \cite{le2024universal} for additional discussion on critical Wasserstein radii in DRO.
For each data provider, we compute the critical $\epsilon_f$ via a line-search process that alters $\epsilon_f$ for each $f$ while keeping all $\epsilon_f',\ f'\neq f$ constant at $0.01$.
Table~\ref{tab:Gridsearch} shows the resulting values.
Notably, the critical $\epsilon_f$ for Cluster 2 is significantly higher than the other clusters. 
We explain this observation with the high electrical distance of the nodes in Cluster 2 from the substation, amplifies the
effects of load and generation on voltage. 
We highlight, that such an analysis can be performed by the DSO to signal a minimum required data quality to data providers or decide on whether datasets are worth obtaining.

\subsubsection{Nodal data provision}
\begin{figure}
    \centering
    \includegraphics[width=0.9\linewidth]{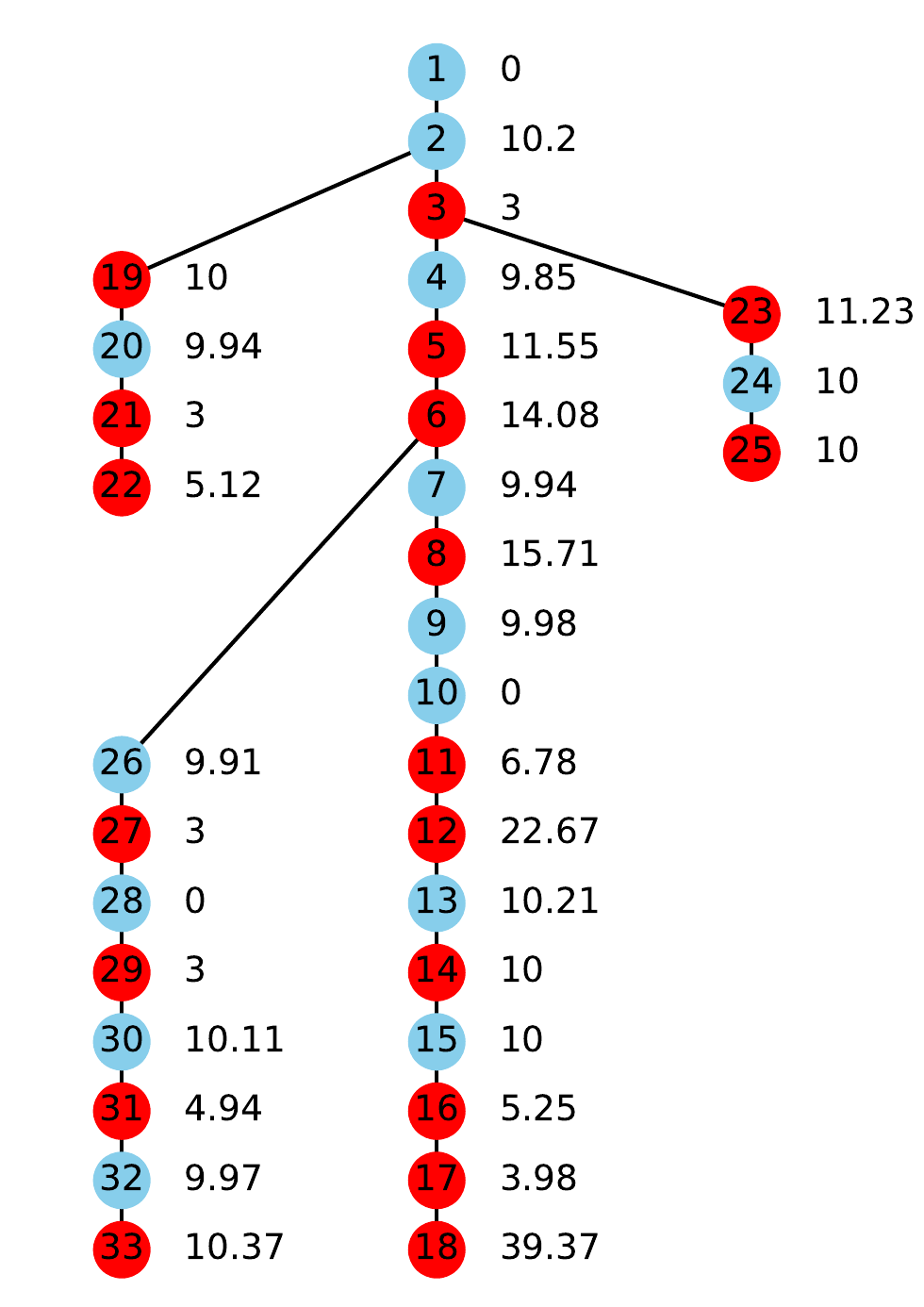}
    \caption{Schematic of the IEEE 33-bus system showing the values of $\mu_f$ for each bus at $\epsilon_f$ $\forall f=1,...,33$ equal 0.01 when the number of clusters is equal to the number of buses in the \textit{High PV}  and \textit{High Load}
    case at 6pm. Nodes marked in red indicate locations with PV.}
    \label{fig:nodal_data_value}
\end{figure}

Finally, Fig.~\ref{fig:nodal_data_value} shows the application of our method for a case where data for each node is supplied by an individual data provider. These results corroborate our previous discussions by showing the higher marginal value of data quality at nodes with higher PV capacity and increased distance from the substation, highlighting the value of high-quality data from the grid edge in the presence of uncertain resources.

\begin{table}[t]
\caption{Critical $\epsilon_f$, $f = 1,...,5$.} 
\label{tab:Gridsearch}
\small
\renewcommand{\arraystretch}{1.2}
\setlength{\tabcolsep}{4pt} 
\centering
\begin{tabular}{|>{\columncolor{lightgray!30}}c!{\vrule width 1.5pt}c|c|c|c|c|}
\hline
\multirow{2}{*} & \multicolumn{5}{c|}{\cellcolor{lightgray!30}{\textbf{f}}} \\
\cline{2-6}
\rowcolor{lightgray!30} 
&\textbf{1} & \textbf{2} & \textbf{3} & \textbf{4} &\textbf{ 5} \\
\noalign{\hrule height 1.5pt}
\textbf{$\bm\lambda^{vol}$}&0.1
 &1
 &0.1
 & 0.005
&0.05
\\
\textbf{$\bm\lambda^{inv}$}& 0.5
& 2.5
& 0.01
& 0.01
& 0.75
\\
\hline
\end{tabular}
\end{table}

\section{Conclusion}
\label{sec:Conclusion}
Motivated by the increasing significance of data-driven decision-making in the power sector, this paper developed a data-driven AC optimal power flow model for radial distribution systems that (i) can accommodate data from various sources and internalize information on data quality and (ii) provides insights on data value via the marginal value of data quality in the final decision.
Taking the perspective of an operator of a distribution system with controllable distributed resources and uncertain load and photovoltaic injections, we first formulated the operator's decision-making problem as a stochastic optimization problem. We then derived a data-driven version of this problem that utilizes data from different data providers (e.g., resource aggregators) and internalizes information on the data quality individual to each data provider. This extends previous results from \cite{mieth2023data}, by allowing data on multiple uncertain parameters to be provided by a single data provider. We then derive a tractable reformulation of the problem. 
Our case study on the IEEE 33-bus test system demonstrates the application of the method and discusses the relationship between data value and data quality.

Data obfuscation for privacy protection, e.g., through differential privacy, is emerging as a prevalent method for data sharing. With our method the data provider can (i) internalize data obfuscation explicitly into the decision-making problem, (ii) compute the maximum allowable obfuscation as the critical $\epsilon$, and (iii) perform an ex-post analysis on the usefulness of a dataset as a function of its quality \textit{without} the necessity of re-running the model.
The decision-dependency of the data value enables the data user to understand in which context data of a certain quality is more (or less) useful.
This offers valuable information for investments in data infrastructure and, at the same time, allows the DSO to transparently communicate data needs as a function of operational circumstance, e.g., to comply with data minimization requirements \cite{ganesh2024data}.

The proposed method also offers potential for data pricing and for applications in machine learning and system state estimation. These applications offer exciting pathways for future research.

\appendices
\section{}
\label{appendix:a}
For given $\bm{\alpha}$, $\bm{q}_c$, $\bm{p}_B$, and $\bm{q}_B$ we can write $g_{\bm{\rho}}(\bm{\alpha},\bm{q}_c,\bm{p}_B,\bm{q}_B,\bm{\delta}) = g_{\bm{\rho}}(\bm{\delta}) = \bm{A}\bm{\delta} + \bm{w}^{vol}$ where $ \bm{A} = \left[ \bm{R}(\bm{I}-{\rm diag}(\bm{\alpha}))\quad  - \bm{R} \quad -\bm{B} \right]\in\mathbb{R}^{N\times3N}$ and $\bm{w}^{vol} = - \bm{R}\bm{p}_B + \bm{B}(\bm{q}_c +\bm{q}_B) +\bm{a}\in\mathbb{R}^{N\times1}$.
Further defining $K=2N$ and $\bm{a_k}$, $c_k$ as the $k$-th row and $k$-th entry of the $K\times3N$ matrix and $K\times1$ vector
\begin{equation*}
    \begin{bmatrix}
        \bm{A} \\ -\bm{A}
    \end{bmatrix}
    \text{ and }
    \begin{bmatrix}
        \bm{w}^{vol} - V_{max} \\ -\bm{w}^{vol} + V_{min} 
    \end{bmatrix},
\end{equation*}
the joint chance constraint \cref{eq:jointccfirst} can be written as 
\renewcommand{\theequation}{A.\arabic{equation}}
\setcounter{equation}{0}
\begin{equation}
    \mathbb{Q}\big\{  \max _{k = 1,...,K}[\langle \bm{a}_k, \bm{\delta}\rangle + c_k ]\le 0 \big\} \ge 1-\eta^{vol}.
\label{eq:compact_joint_cc}
\end{equation}
Using samples $\{\widehat{\bm{\delta}}_i\}_{i=1}^{I}$ Eq.~\cref{eq:compact_joint_cc} can be solved using the data-driven approach from \cite{dall2017chance}, which leads to the result in \cref{eq:joint_cc_sample}. 

\section{Chance constraint reformulation}
\label{appendix:b}
\subsection{Derivation of $\eqref{eq:maincvar0}$}

We first reformulate \cref{eq:voltagecvar1} as\cite{rockafellar2000optimization}:
\begin{align}
\begin{cases}
    \varphi^{vol} \le 0  \\    \varpi^{vol} +\varphi^{vol} \le 0 \\   \eta^{vol} \varpi^{vol} \ge \sup_{\mathbb{Q}\in \mathcal{A}} \mathbb{E}_\mathbb{Q}[ \max _{k = 1,...,K+1}\langle \bm{a'}_k, \bm{\delta}\rangle + c'_k],  
\end{cases}\label{eq:jointcvar11}
\end{align}
where $\bm{a'}_k= \bm{a}_k,$ $k = 1,...,K$, $\bm{a'}_{K+1} = \bm{0}_D$, $c'_k= c_k - \varphi^{vol}$ $k = 1,...,K$, $c'_{K+1} = 0$ such that $\bm{0}_F$ denotes a vector of zeros of length $F$. Variable $\varpi^{vol}$ is an auxiliary decision variable. We then use \cref{eq:mwsdro_general_formulation_stand_data} to reformulate the worst-case expectation appearing in \cref{eq:jointcvar11} as
\begin{flalign}
&\sum_{f=1}^{F} \lambda_{f}^{vol} \epsilon_{f} + \frac{1}{I}\sum_{i=1}^{I} s_{i}^{vol} \le \eta^{vol} \varpi^{vol} \\
&s_{i}^{vol}\! \ge \!\sup_{\bm{\delta}\in \Xi} \big[ \!\max _{k = 1,...,K+1}(\langle \bm{a}^\prime_k, \bm{\delta}\rangle + c^\prime_k)\! -\!\!\sum_{f=1}^{F} \!\!\lambda_{f}^{vol}\|\bm{\delta}_{f} \!-\! \hat{\bm{\delta}}_{{f},{i}}\|\big] \nonumber\\& \qquad  i=1,\ldots,I \label{eq:jointcvar4}
\end{flalign}
We have used $p=1$ in this reformulation. 
Further, using $\lambda_{f}^{vol}\|\bm{\delta}_f - \widehat{\bm{\delta}}_{f,i}\|=\sum_{m=1}^{3N_f} \lambda_{f}^{vol}|\delta_{f,m}-\widehat{\delta}_{f,i,m}|$ and introducing $z_{k,f,i}$ as an auxiliary variable we can rewrite \eqref{eq:jointcvar4} as 
\begin{flalign}
&  s_{i}^{vol}\!\!\ge\!\!  c^\prime_k \!\!+\!\! \!\!\!\!\!\!\!\!\!\!\sup_{{\scriptstyle \delta_{f,m}\in[\underline{\delta}_{f,m},\overline{\delta}_{f\!,m}]}}\! \sum_{f=1}^{F}  \!\sum_{m=1}^{3N_{f}} \!(\!z_{k,f,i} \widehat{\delta}_{f,i,m} \!+\! (\!a'_{k,f,m}\!\! \! -\!z_{k,f,i}) \delta_{f,m}\!) \label{eq:reform_10} \\ 
&k=1,...,K+1   \quad i=1,\ldots,I\\
& |z_{k,f,i}| \le\lambda_{f}^{vol} \quad k=1,...,K+1   \quad i=1,\ldots,I
\end{flalign}
Finally, we resolve the inner supremum in \cref{eq:reform_10} using 
\begin{align}
    & \sup_{\delta_{f,m}\in[\underline{\delta}_{f,m},\overline{\delta}_{f,m}]}\big((a'_{k,f,m} - z_{k,f,i} ) \delta_{f,m} \big) \\
=   & \!\!\!\inf_{\substack{ u_{k,f,i,m},  l_{k,f,i,m}\ge 0:\\(a'_{k,f,m} - z_{k,f,i}) = (u_{k,f,i,m}- l_{k,f,i,m})}}\!\!\! u_{k,f,i,m} \overline{\delta}_{f,m} - l_{k,f,i,m} \underline{\delta}_{f,m},\nonumber
\end{align}
which result in $\eqref{eq:maincvar0}$.

\subsection{Derivation of $\eqref{eq:maincvarinv}$}
We reformulate \cref{eq:thirddcvar} similarly to \cref{eq:voltagecvar1} by introducing auxiliary variable $\varpi_n^{inv}$ and then writing
\begin{align}
\begin{cases}
      \varphi_n^{inv} \le 0  \qquad n = 1,...,N  \\
     \varpi_n^{inv} +  \varphi_n^{inv}  \le 0  \qquad  n = 1,...,N  \\
     \eta^{inv} \varpi_n^{inv}\ge \!\sup_{\mathbb{Q}\in \mathcal{A}} \mathbb{E}_\mathbb{Q}[\max \{((1-\!\alpha_n) p_{av,n})^2 +\! w_{n}^{inv}]^+ \\ \qquad  n = 1,...,N. 
\end{cases}\!\!\label{eq:thirrd-lastline}
\end{align}
where $w_{n}^{inv} = q_{c,n}^2 - S_n^{2} - \varphi_n^{inv} , n=1,...,N$. Using \cite[Propostion~1]{mieth2023data}, we write objective \cref{eq:msw_dro_objective} as:
Applying \cite[Proposition 2]{mieth2023data} to the worst-case expectation in \eqref{eq:thirrd-lastline} then leads to \cref{eq:maincvarinv}.

\section{}
\label{appendix:c}
 \allowdisplaybreaks
\begin{subequations}
\begin{flalign}
&\min  \sum_{f=1}^{F} \lambda_f^{co} \epsilon_f +  \sum_{n=1}^{N}( (\frac{1}{I} \sum_{i=1}^{I} s_{n,i}^{co1}) +  e_n ||q_{c,n}|+|q_{B,n}||\nonumber\\& +(\frac{1}{I} \sum_{i=1}^{I} s_{n,i}^{co2})  )\\
&\text{s.t.}  \quad n=1,\dots,N \quad i=1,\dots,I: \nonumber \\
  &s_{n,i}^{co1}\! \ge c_n ((\bm r_n \overline{\bm \delta}_{f(n)})\! -\! (1-\alpha_n) (\bm m_n \underline{\bm \delta}_{f(n)}) - p_{B}) - \lambda_{f(n)}\nonumber\\&  (\!((\bm r_n \overline{\bm \delta}_{f(n)})\! - (\bm r_n \hat{\bm \delta}_{f(n),i})) -((\bm m_n \underline{\bm \delta}_{f(n),i}) -\! (\bm m_n \hat{\bm \delta}_{f(n),i}))\!) \label{objfinal1}\\
 &s_{n,i}^{co1} \ge d_n ((1-\alpha_n) (\bm m_n\overline{\bm \delta}_{f(n)}) - (\bm r_n \underline{\bm \delta}_{f(n)}) + p_{B}) - \lambda_{f(n)}\nonumber\\ &  ((- (\bm r_n \underline{\bm \delta}_{f(n)}) \!+ (\bm r_n \hat{\bm \delta}_{f(n),i})) +\! ((\bm m_n\overline{\bm \delta}_{f(n)}\!) -\! (\bm m_n \hat{\bm \delta}_{f(n),i}\!))\!)\\
 & s_{n,i}^{co1} \ge c_n ((\bm r_n \hat{\bm \delta}_{f(n),i}) - (1-\alpha_n) (\bm m_n \hat{\bm \delta}_{f(n),i}) - p_{B}) \\
& s_{n,i}^{co1} \ge  d_n ((1-\alpha_n) (\bm m_n\hat{\bm \delta}_{f(n),i}) - (\bm r_n \hat{\bm \delta}_{f(n),i}) + p_{B})\\
& s_{n,i}^{co1} \ge 0\\
 &s_{n,i}^{co2} \!\ge h_n \alpha_n (\bm m_n\overline{\bm \delta}_{f(n)})  - \lambda_{f(n)} ((\bm m_n\overline{\bm \delta}_{f(n)}) -(\bm m_n \hat{\bm \delta}_{f(n),i})) \\
&s_{n,i}^{co2}\! \ge h_n \alpha_n (\bm m_n\underline{\bm \delta}_{f(n)} ) + \lambda_{f(n)} ((\bm m_n\underline{\bm \delta}_{f(n)}) - (\bm m_n \hat{\bm \delta}_{f(n),i}))\\
 &s_{n,i}^{co2} \!\ge h_n \alpha_n (\bm m_n \hat{\bm \delta}_{f(n),i}) \\
& s_{n,i}^{co2} \ge 0. \label{objfinal2}
\end{flalign}%
\end{subequations}%
\allowdisplaybreaks[0]%
Here, $\bm r_n$ and $\bm m_n$ are column vectors with the same dimensions as $\bm{\delta}_{f(n)}$. 
Vector $\bm r_n$ has all zero entries except at the index matching $p_{l,n}$ in $\bm{\delta}_{f(n)}$.
Similarly, vector $\bm m_n$ has all zero entries except at the index matching $p_{av,n}$ in $\bm{\delta}_{f(n)}$.

\bibliographystyle{IEEEtran}
\bibliography{literature}

\end{document}